\documentclass[pdflatex,sn-basic]{sn-jnl}

\usepackage[T1]{fontenc}
\usepackage[utf8]{inputenc}
\usepackage{lmodern}
\usepackage{indentfirst}
\input{glyphtounicode}
\usepackage{graphicx}
\graphicspath{{figures/}}
\usepackage{amsmath,amssymb,amsfonts}
\usepackage{booktabs}
\usepackage{tabularx}
\usepackage{array}
\usepackage{placeins}
\usepackage{longtable}
\usepackage{needspace}
\usepackage{etoolbox}
\hypersetup{bookmarksdepth=2}

\AtBeginEnvironment{thebibliography}{\interlinepenalty=10000}
\usepackage{algorithm}
\usepackage{algpseudocode}
\usepackage{xr}
\hypersetup{pdfauthor={Jue Wang, Ximo Li, Xiaofei Xu},hypertexnames=false}
\newcommand{\mat}[1]{\boldsymbol{#1}}
\providecommand{\doiurl}[1]{\url{https://doi.org/#1}}

\begin{document}

\title[A Sharp Threshold for Exact Label Recovery]{A Sharp Threshold for Exact Label Recovery in Incomplete Binary Latent Class Models}

\author[1]{\fnm{Jue} \sur{Wang}\enspace
  \href{https://orcid.org/0009-0007-6618-5190}{\textnormal{\scriptsize\mbox{(0009-0007-6618-5190)}}}}
\email{jue.wang@whu.edu.cn}
\author[2]{\fnm{Ximo} \sur{Li}\enspace
  \href{https://orcid.org/0009-0005-6458-5890}{\textnormal{\scriptsize\mbox{(0009-0005-6458-5890)}}}}
\email{xmli007@stu.xidian.edu.cn}
\author*[1]{\fnm{Xiaofei} \sur{Xu}\enspace
  \href{https://orcid.org/0000-0002-8019-9245}{\textnormal{\scriptsize\mbox{(0000-0002-8019-9245)}}}}
\email{xiaofeix@whu.edu.cn}
\affil[1]{\orgdiv{School of Mathematics and Statistics}, \orgname{Wuhan University},
  \orgaddress{\city{Wuhan}, \state{Hubei}, \country{China}}}
\affil[2]{\orgdiv{School of Economics and Management}, \orgname{Xidian University},
  \orgaddress{\city{Xi'an}, \state{Shaanxi}, \country{China}}}

\abstract{Latent class models help identify groups with distinct response patterns in surveys and educational assessments.
Existing theory for spectral clustering with likelihood refinement focuses on complete responses, leaving open how much observed information is needed to recover every label when class response probabilities must be estimated from incomplete data.
We construct and analyze a cross-fitted estimator that uses spectral clustering to estimate class-specific response probabilities in the training data and assigns labels to held-out subjects using the observed-data likelihood.
Under independent, item-specific missingness completely at random, we measure recoverability by the minimum pairwise Chernoff information between the class-specific distributions of a subject's observed responses.
We prove that, in a regular high-dimensional regime with comparable observation probabilities, the estimator recovers the true partition with probability tending to one when this information exceeds $\log n$, the logarithm of the sample size, by a fixed relative margin.
We establish a matching lower bound showing that, below $\log n$ by a fixed relative margin, the minimax probability of exact recovery tends to zero even with known model parameters and class sizes.
In a paired simulation with equal class priors and information at 1.3 times the threshold, replacing the true response probabilities with estimates from the true training groups reduced exact recovery from 99 to 52 out of 100 datasets.
Learning the response probabilities can therefore preserve the oracle's first-order threshold while still affecting exact recovery at finite sample sizes.}
\keywords{Latent class analysis; Model-based clustering; Spectral clustering; Missing data; Cluster recovery; High-dimensional binary data}
\pacs[MSC Classification]{62H30 (primary); 62H25; 62H12}

\maketitle

\section{Introduction}

Surveys and educational assessments often use responses to multiple items to
identify groups with distinct response patterns. When response records are
incomplete, the information available for recovering these groups also depends
on which items are observed. The record for a binary item distinguishes a missing
response, an observed zero, and an observed one. Filling missing entries with
zeros and discarding the observation mask merges the first two outcomes and can
reduce the information available for distinguishing classes. We study binary
latent class models under independent, item-specific missingness completely at
random (MCAR) and ask how much observed information is needed to recover the
true partition when class response probabilities must be estimated from
incomplete data.

Latent class models represent multivariate responses as finite mixtures of
product distributions \citep{Goodman1974}. \citet{Allman2009} establish generic
identifiability results for such models. Under ignorable missingness, the
observed-data likelihood sums the probabilities of all complete response
patterns consistent with the observed record \citep{Rubin1976,LittleRubin2019}.
The expectation--maximization (EM) algorithm \citep{DempsterEtAl1977} is a
standard approach to fitting these models and is used in the
\texttt{poLCA} package \citep{LinzerLewis2011}. Initialization can affect
likelihood fitting in mixture models \citep{BaudryCeleux2015}. For binary
responses, \citet{NajafiEtAl2020} study clusterability in fully observed
Bernoulli mixtures, and \citet{ZengGuXu2023} use a tensor method to initialize
EM for large-scale latent class analysis.

Spectral methods provide another way to obtain initial groups.
\citet{ZhuWangSamworth2022} study high-dimensional subspace estimation under
heterogeneous missingness, and \citet{Qing2026} gives guarantees for
regularized spectral clustering of ordered polytomous responses. Low-rank
center-separation arguments connect spectral approximation error to
classification error \citep{LofflerZhangZhou2021,ZhangZhou2024}.
\citet{LeiZhu2017} use sample splitting to refine initial groups in network
clustering. For fully observed binary responses, \citet{LyuGu2026} combine
spectral initialization with one likelihood update and analyze a version with
sample splitting and label alignment. They establish a minimax exponential
rate for expected classification loss and give sufficient conditions for
exact recovery.

Applying this approach to incomplete responses raises the question of whether
response probabilities estimated from the initial groups preserve the
information that distinguishes classes. We use cross-fitting to form spectral
groups and estimate their response probabilities on each training set, then use
those estimates and the observed-data likelihood to classify held-out subjects.
The training partition determines whose responses enter each profile estimate,
so the profile error depends on how the same training subjects are grouped.
We show that the combined effect of profile estimation errors across items
reduces the classification error exponent by only a lower-order term,
preserving the leading exponent obtained with known class profiles.

The hardest class pair determines this exponent. Let $C_*^{\mathrm{obs}}$
denote the minimum pairwise Chernoff information between the class-specific
distributions of one subject's observed responses. With a fixed number of
classes whose proportions are bounded away from zero, item counts of the same
order as the sample size $n$, comparable observation probabilities, and
regularly scaled and separated profiles, the sharp first-order threshold for
exact recovery is $C_*^{\mathrm{obs}}=\log n$. Above the threshold by a fixed
relative margin, both centered mean imputation and propensity-weighted
residuals recover the true partition with probability tending to one when
followed by the same likelihood refinement. Below the threshold by a fixed
relative margin, the minimax probability of exact recovery tends to zero even
with known class profiles, observation probabilities, and class sizes. For
more heterogeneous observation probabilities, we derive finite-sample
recovery bounds and extend the analysis to settings with rarely observed
items or random retention near the screening threshold. The simulations
reveal how recovery can differ at finite sample sizes under a common
first-order threshold. In the asymmetric design with observed information at
1.3 times the threshold, replacing the true profiles with estimates from the
true training groups reduced the number of exact recoveries from 99 to 52 out
of 100 datasets. This shows that estimating response probabilities can
substantially affect finite-sample recovery even when the training groups are
correct. For students in held-out schools in the US and Singapore PISA 2022
samples, latent class models fitted from either spectral start achieved higher
response prediction scores than Rasch-SD.

Section~\ref{sec:method} describes the model, observed information, and
estimator. Section~\ref{sec:recovery-theory} gives finite-sample recovery bounds
(Theorems~1 and~2) and the sharp first-order threshold (Theorem~3), and
Section~\ref{sec:extensions} develops the informative-core and random-retention
extensions. Sections~\ref{sec:simulation}, \ref{sec:pisa}, and
\ref{sec:discussion} present the simulations, the PISA application, and the
discussion. Appendices~\ref{app:observed-information}--\ref{app:information}
prove the information propositions,
Appendices~\ref{app:spectral}--\ref{app:sharp-outline} give the recovery proofs
and key estimates, and
Appendices~\ref{app:common-threshold-proof}--\ref{app:random-retention-proof}
prove the recovery corollaries. Online Resource~1 contains longer technical
arguments and the encoding and screening extensions
(Sections~\ref*{app:sharp-proof}--\ref*{app:encoding-extensions}), as well as
larger-scale experiments and additional diagnostics
(Sections~\ref*{app:diagnostics}--\ref*{app:sharp-experiment}).

\section{Model, information, and estimator}
\label{sec:method}

\subsection{Item-specific MCAR latent class model}
\label{sec:model-assumptions}

For $n$ subjects, $J$ items, and a fixed number $K$ of classes, let
$\boldsymbol z=(z_i)_{i=1}^n$ be fixed labels with
$z_i\in\{1,\ldots,K\}$. Conditional on $z_i=k$, the response $X_{ij}$ has
class-specific response probability $\theta_{kj}$, while the observation
indicator $\Omega_{ij}$ has item-specific observation probability $\pi_j$.
Conditional on $\boldsymbol z$, all response and mask cells are mutually
independent. Writing $\operatorname{Ber}(p)$ for a Bernoulli law with success
probability $p$,
\begin{align*}
X_{ij}\mid z_i=k&\sim\operatorname{Ber}(\theta_{kj}),&
\Omega_{ij}&\sim\operatorname{Ber}(\pi_j),&
\mat{\Omega}&\perp(\mat{X},\boldsymbol{z}).
\end{align*}

The observed symbol $W_{ij}=(\Omega_{ij},\Omega_{ij}X_{ij})$ records a missing
response, an observed zero, or an observed one as $(0,0),(1,0),(1,1)$,
respectively. In class $k$, these three outcomes have probabilities
$1-\pi_j$, $\pi_j(1-\theta_{kj})$, and $\pi_j\theta_{kj}$. Let $P_k$ denote
the product law of the observed row $\boldsymbol W_i=(W_{ij})_{j=1}^J$;
$E_k$ and $\operatorname{Var}_k$ denote expectation and variance under this
law. The profile matrix is $\mat\Theta=(\theta_{kj})$, with rows
$\boldsymbol\theta_k$. The class sizes
$n_k=\sum_i\mathbf1\{z_i=k\}$ and proportions $\alpha_k=n_k/n$ form the
vector $\boldsymbol\alpha$ and satisfy
$0<\alpha_{\min}\leq\alpha_k\leq\alpha_{\max}<1$.

Under regular probability scaling, write
$\theta_{kj}=\rho q_{kj}$, where $\rho>0$ and the scaled coefficients satisfy
\(0<q_-\leq q_{kj}\leq q_+<\infty\) and
\(\rho q_+\leq1-\epsilon_0\). The common scale $\rho$ may tend to zero or
remain bounded away from zero. Here $q_-,q_+,\epsilon_0>0$ are fixed. Only
the sparse expansion in Appendix~\ref{app:observed-information} requires $\rho\to0$.

For cross-fitting, use a data-independent uniform partition into $L$ target
folds $F_1,\ldots,F_L$ with pre-fixed, nearly equal sizes. The training set
for fold $f$ is $T_f=\{1,\ldots,n\}\setminus F_f$; write
$m_f=|T_f|$ and $m=\min_fm_f$. All probability statements condition on the
fixed label vector and parameter sequence. Randomness includes responses,
observation masks, the random fold assignment, and algorithmic seeding.

\subsection{Observed-data information}
\label{sec:chernoff-information}

An item distinguishes classes $k$ and $\ell$ only when its response is
observed. Following \citet{VanErvenHarremoes2014}, write $b_{k\ell j}$ for
the complete-data affinity and $a^{\mathrm{obs}}_{k\ell j}$ for the affinity
after accounting for the common missingness probability $1-\pi_j$:
\begin{align*}
b_{k\ell j}
&=\sqrt{\theta_{kj}\theta_{\ell j}}
+\sqrt{(1-\theta_{kj})(1-\theta_{\ell j})},&
a^{\mathrm{obs}}_{k\ell j}
&=1-\pi_j+\pi_jb_{k\ell j}.
\end{align*}

\textbf{Proposition 1 (exact observed R\'enyi information).} Since the row
distributions are products, their order-one-half R\'enyi divergence is the
sum of the item divergences. Taking the logarithm of each observed-data
affinity with the convention $-\log 0=+\infty$ gives
\begin{equation}
I^{\mathrm{obs}}_{k\ell}=\sum_{j=1}^{J}-2\log a^{\mathrm{obs}}_{k\ell j}.
\label{eq:information}
\end{equation}

The Chernoff information uses a common order $0<s<1$ across items. Defining
each coefficient by continuous extension in the response probabilities at
fixed $s$, we obtain
\begin{align*}
c_{k\ell j}(s)
&=(1-\pi_j)+\pi_j\bigl\{\theta_{kj}^{s}\theta_{\ell j}^{1-s}
+(1-\theta_{kj})^{s}(1-\theta_{\ell j})^{1-s}\bigr\},\\
C^{\mathrm{obs}}_{k\ell}
&=-\inf_{0<s<1}\sum_{j=1}^{J}\log c_{k\ell j}(s).
\end{align*}

At $s=1/2$ the coefficient equals the observed-data affinity, so
$C^{\mathrm{obs}}_{k\ell}\ge I^{\mathrm{obs}}_{k\ell}/2$. Under regular scaling
the coefficient extends to both endpoints with value one. At boundary
response probabilities we take the infimum over the open interval $(0,1)$.

For a subset $S$ of the $J$ items, write $I^{\mathrm{obs}}_{k\ell,S}$ and
$C^{\mathrm{obs}}_{k\ell,S}$ for the information in the corresponding product
distributions, obtained by restricting each sum to $j\in S$. A star in place
of the class indices denotes the minimum over distinct pairs. Thus
$C^{\rm obs}_{*,S}$ measures how well the items in $S$ distinguish the hardest
pair under the Chernoff criterion, and $I^{\rm obs}_{*,S}$ is the analogous
minimum for R\'enyi divergence. Omitting $S$ includes all $J$ items.

\paragraph{Zero filling.}

The observed-data likelihood uses the full symbol
$\mat{W}=(\mat{\Omega},\mat{\Omega}\odot\mat{X})$, where $\odot$ denotes
entrywise multiplication. The zero-filled input,
in contrast, keeps only $\mat{V}=\mat{\Omega}\odot\mat{X}$, so that the mask
is not retained. Its entries satisfy
$V_{ij}\mid z_i=k\sim\operatorname{Ber}(\pi_j\theta_{kj})$.
Write $I^+_{k\ell,S}$ for the order-one-half R\'enyi divergence between
these zero-filled product distributions on $S$.

\textbf{Proposition 2 (information loss under zero filling).}
Zero filling cannot increase the order-one-half information:
$I^+_{k\ell,S}\le I^{\mathrm{obs}}_{k\ell,S}$. Itemwise equality holds
exactly when $\theta_{kj}=\theta_{\ell j}$, $\pi_j=0$, or $\pi_j=1$.
Under regular scaling, for sufficiently small $\rho$,
\[
0\le I^{\mathrm{obs}}_{k\ell,S}-I^+_{k\ell,S}
\le C\rho I^{\mathrm{obs}}_{k\ell,S}.
\]

Appendix~\ref{app:observed-information} proves Propositions~1 and~2 and
gives the exact information-loss identity and its sparse expansion.
The recovery benchmark for the class pair
$k,\ell$ is $C^{\mathrm{obs}}_{k\ell}$.

\paragraph{Regular probability scaling.}
The observation-weighted profile contrast is
\begin{equation*}
Q_{k\ell}
=\sum_{j=1}^{J}\pi_j
\frac{(\theta_{kj}-\theta_{\ell j})^2}{\rho}
=\rho\sum_{j=1}^{J}\pi_j(q_{kj}-q_{\ell j})^2.
\end{equation*}

\textbf{Proposition 3 (regular-scaling information equivalence).} Under the
regular probability bounds above, there are fixed positive constants $c$ and
$C$ such that
$cQ_{k\ell}\leq I^{\mathrm{obs}}_{k\ell}\leq CQ_{k\ell}$ and
$C^{\mathrm{obs}}_{k\ell}\asymp Q_{k\ell}$. Here $a\asymp b$ means
$cb\le a\le Cb$ for fixed constants $c,C>0$.
Appendix~\ref{app:information} contains the proof.

\subsection{Estimated preprocessing and spectral partition}

For a training row set $T$, restrict all matrices in this subsection to those
rows. Given a pre-specified threshold $\tau\in[0,1]$, estimate the observation
probabilities and item centers, then screen the items:
\[
\widehat\pi_j(T)=|T|^{-1}\sum_{i\in T}\Omega_{ij},\qquad
S(T)=\{j:\widehat\pi_j(T)\ge\tau,\widehat\pi_j(T)>0\},
\]
\[
\widehat c_j=\frac{\sum_{i\in T}\Omega_{ij}X_{ij}+1/2}
{\sum_{i\in T}\Omega_{ij}+1}.
\]

Here $\widehat c_j$ is the Jeffreys-smoothed item center. We write
$\widehat\pi_j^{(f)}=\widehat\pi_j(T_f)$ and suppress the training-set index
on nuisance estimates when it is unambiguous; $\tau$ may vary with sample
size. The exponent $\gamma\in[0,1/2]$ determines the propensity weighting in
\begin{equation}
\widehat{\mat{R}}_{\gamma,T}=
\mat{H}_T\left[\mat{\Omega}\odot(\mat{X}-\mathbf1\widehat{\boldsymbol c}^\top)\right]_{S(T)}
\operatorname{diag}(\widehat\pi_{S(T)}^{-\gamma})
.
\label{eq:residual}
\end{equation}

Here $\mathbf1$ is the all-ones vector of length $|T|$ and
$\mat{H}_T=\mat{I}_{|T|}-\mathbf1\mathbf1^\top/|T|$ centers the training rows.
The subscript $S(T)$ restricts a matrix to those columns or a vector to
those entries.

The endpoint $\gamma=0$ corresponds to centered mean imputation (Mean), whereas
$\gamma=1/2$ corresponds to propensity-weighted residuals (Residual).
Filling missing entries with $\widehat c_j$ and then centering and weighting
the columns gives the same matrix as Equation~\eqref{eq:residual}.

For $T=T_f$, set $\widehat S_f=S(T_f)$, abbreviate this set by
$\widehat S$ when the fold is clear, and define
$r_f=\min\{K-1,|\widehat S_f|,m_f-1\}$. Compute the exact singular value
decomposition (SVD), retaining the leading left singular vectors
$\widehat{\mat U}_{r_f}$ and the diagonal matrix of singular values
$\widehat{\mat\Sigma}_{r_f}$. Their product defines the spectral score rows;
when the input and rank are clear, write this product as $\mat U\mat\Sigma$.
Use $R\ge1$ independent D-squared starts: choose the first center uniformly,
then sample subsequent centers with probability proportional to the squared
distance from the nearest selected center \citep{ArthurVassilvitskii2007}.
Lloyd updates alternate nearest-center assignment and centroid updates.
Retain the solution with minimum inertia, the within-cluster sum of squared
distances. Keeping the singular-value scales links low-rank reconstruction
error to class-center separation. Algorithm~\ref{alg:crossfit}
specifies the full procedure. Appendix~\ref{app:exact-conventions} gives the
exact-arithmetic conventions.

\subsection{Profile estimation and likelihood refinement}

Let $\widehat z_i^{(f)}$ be the spectral training label for subject $i\in T_f$.
Define the estimated membership indicator by
$\widehat Z_{ik}^{(f)}=\mathbf1\{\widehat z_i^{(f)}=k\}$. The Jeffreys-smoothed
profile and class-proportion estimates are
\begin{equation}
\widehat\theta^{(f)}_{kj}
=\frac{\frac12+\sum_{i\in T_f}\widehat Z_{ik}^{(f)}
\Omega_{ij}X_{ij}}
{1+\sum_{i\in T_f}\widehat Z_{ik}^{(f)}\Omega_{ij}},
\qquad
\widehat\alpha^{(f)}_k
=\frac{\frac12+\sum_{i\in T_f}\widehat Z_{ik}^{(f)}}
{m_f+K/2}.
\label{eq:profiles}
\end{equation}

The true training proportion is
$\alpha_{k,f}=m_f^{-1}\sum_{i\in T_f}\mathbf1\{z_i=k\}$; it appears only
in the analysis. For each held-out subject, the classifier scores each class
by summing the Bernoulli log-likelihood terms for that subject's observed
responses on the retained items in $\widehat S_f$ and adding the logarithm of
the estimated class proportion. It assigns the class with the highest score.
Factors depending only on the mask cancel in this comparison,
but the mask still determines which response terms enter the score.
For probabilities $p,q$, use the unhalved squared Hellinger distance
$d_{\rm Ber}^2(p,q)=(\sqrt p-\sqrt q)^2+(\sqrt{1-p}-\sqrt{1-q})^2$.
For profiles $\boldsymbol u,\boldsymbol v$ on any item set $S$, define
\[
d_S^2(\boldsymbol u,\boldsymbol v)=\sum_{j\in S}d_{\rm Ber}^2(u_j,v_j).
\]

The Hungarian algorithm aligns the class labels across folds by minimizing
the sum of squared Hellinger distances. Concatenating the resulting target predictions
defines \textsc{ResidualCrossfit}; setting $\gamma=0$ yields
\textsc{MeanCrossfit}.

\begin{algorithm}[!htbp]
\normalsize
\caption{Residual spectral initialization with cross-fitted likelihood refinement}
\label{alg:crossfit}
\begin{algorithmic}[1]
\Require Observed data $(\mat{\Omega},\mat{\Omega}\odot\mat{X})$, fixed class count $K\ge2$ and fold count $2\le L\le n$, threshold $\tau\in[0,1]$, restart count $R\ge1$, exponent $\gamma\in[0,1/2]$, recorded realizations of data-independent random seeds.
\State Partition the subjects into nearly equal target folds $F_1,\ldots,F_L$, independently of the data.
\For{$f=1,\ldots,L$}
\State Set $T_f=\{1,\ldots,n\}\setminus F_f$; estimate $\widehat\pi_j$ and Jeffreys item centers $\widehat c_j$ using only $T_f$.
\State Return $\bot$ (failure) if $K>|T_f|$. Screen the items and retain $\widehat S_f=\{j:\widehat\pi_j\ge\tau,\widehat\pi_j>0\}$; return $\bot$ if $\widehat S_f$ is empty.
\State Form $\widehat{\mat{R}}_{\gamma,f}=\mat{H}_{T_f}[\mat{\Omega}\odot(\mat{X}-\mathbf1\widehat{\boldsymbol c}^{\top})]_{\widehat S_f}\operatorname{diag}(\widehat\pi_{\widehat S_f}^{-\gamma})$ on $T_f$.
\State Set $r_f=\min\{K-1,|\widehat S_f|,|T_f|-1\}$; compute the leading $r_f$ spectral scores $\widehat{\mat{U}}_{r_f}\widehat{\mat{\Sigma}}_{r_f}$.
\State Run $R$ independent D-squared/Lloyd starts; retain the training partition with minimum inertia, using the stated tie and empty-center rules.
\State Estimate Jeffreys class profiles $\widehat{\mat{\Theta}}^{(f)}$ and proportions $\widehat{\boldsymbol\alpha}^{(f)}$ by Equation~\eqref{eq:profiles}.
\EndFor
\State Return $\bot$ if any $\widehat S_f\cap \widehat S_1=\varnothing$ ($f>1$); otherwise align to fold 1 by Hungarian matching of summed squared Hellinger distances on $\widehat S_f\cap \widehat S_1$, applying the permutation to both $\widehat{\mat\Theta}^{(f)}$ and $\widehat{\boldsymbol\alpha}^{(f)}$.
\For{$f=1,\ldots,L$ and $i\in F_f$}
\State Assign the aligned class label maximizing $\log\widehat\alpha_k^{(f)}+\sum_{j\in \widehat S_f:\Omega_{ij}=1}\{X_{ij}\log\widehat\theta_{kj}^{(f)}+(1-X_{ij})\log(1-\widehat\theta_{kj}^{(f)})\}$.
\EndFor
\State \Return Concatenated target labels, up to one global permutation.
\end{algorithmic}
\end{algorithm}
 \FloatBarrier

\section{Recovery theory}
\label{sec:recovery-theory}

The main result is the first-order boundary $C_*^{\rm obs}=\log n$ for
recovering all labels. Under comparable observation probabilities,
Theorem~3 and Corollary~1 show that the fitted procedure meets this boundary
at fixed positive margins and that even an estimator with known model
parameters and class sizes cannot recover uniformly below it.

We analyze the same estimator in two ways. Section~\ref{sec:general-route}
gives finite-sample bounds that allow for itemwise heterogeneity in noise and
separation, using uniform bounds on itemwise profile errors to analyze held-out
classification. Section~\ref{sec:sharp-recovery} uses aggregate profile error
across items to obtain the sharp threshold.

Throughout this section, constants may depend on the fixed class and fold
counts, class balance, profile bounds, and the stated comparison parameters.
Dimensions, probability scales, thresholds, and error tolerances may vary
with $n$. Unless stated otherwise, confidence parameters lie in $(0,1)$.
Algorithmic seeding is independent of responses and masks before execution.

Let $\mathfrak S_K$ be the set of permutations of $\{1,\ldots,K\}$.
For estimated and true label vectors, define normalized permutation loss by
\begin{equation}
\operatorname{err}(\widehat{\boldsymbol z},\boldsymbol{z})
=\frac1n\min_{\sigma\in\mathfrak S_K}
\sum_{i=1}^n\mathbf1\{\sigma(\widehat z_i)\ne z_i\}.
\label{eq:loss}
\end{equation}

This definition permits unused estimated labels. For a failed return, set
$\operatorname{err}(\bot,\boldsymbol z)=1$. We write
$P_{\rm perm}\{\widehat{\boldsymbol z}\ne\boldsymbol{z}\}$ for
$P\{\operatorname{err}>0\}$; its complement is the exact-recovery
probability. The training classification loss $\operatorname{Mis}_\gamma$
uses the same definition with the training-sample size in place of $n$.

We refer to the conditions introduced in Section~\ref{sec:model-assumptions}
as follows: \textbf{A1} denotes fixed balanced classes; \textbf{A2}, mutually
independent item-specific MCAR; \textbf{A3}, regular probability scaling; and
\textbf{A4}, a data-independent partition into a fixed number of nearly equal
folds.

\subsection{General recovery from estimated preprocessing}
\label{sec:general-route}
For nonnegative quantities, $a\lesssim b$ means $a\le Cb$ for a constant
$C>0$ independent of $n$.
The recovery analysis has two stages. First, we bound training
misclassification in terms of preprocessing perturbation and class-center
separation. We then use the errors of the profiles estimated from the training
groups to bound the probability of recovery failure on the held-out rows.
The following display gives the main terms of both bounds.
\begin{align*}
\text{training classification loss}
&\lesssim
\frac{\text{preprocessing perturbation}^{2}}
{\text{training size}\times\text{center separation}^{2}},\\
\Pr(\text{exact-recovery failure})
&\lesssim \Pr(\text{required events fail})\\
&\quad+n\exp\{-\text{observed information}
+\text{estimation penalty}\}.
\end{align*}

Theorem~1 gives the training bound with perturbation scale $b_\gamma$ and
separation $\Delta_{\gamma,c}$. Theorem~2 gives the recovery bound under
profile-matching and concentration conditions, with training loss entering
the profile-error scale $r_m$ and hence the score penalty $R_m^J$. The bound
adds the probability of a failed training or held-out count event to an
exponential bound on misclassification with all events satisfied. The
Chernoff exponent in the latter term is at least the core's observed
Chernoff information minus $2R_m^J$.

\paragraph{Retention and scaling.}

Set the screening threshold to $\tau=h>0$, allowing $h$ to vary with sample
size. Define $S_t=\{j:\pi_j\ge t\}$; then $S_{2h}$ is the population core
and $S_{h/2}$ is the outer set. Write $\widehat S_{h,f}=\widehat S_f$ to
make the threshold explicit, and write $\widehat S_h$ when the fold index is
suppressed.
Threshold subscripts on information abbreviate the corresponding set,
for example $I^{\rm obs}_{*,2h}=I^{\rm obs}_{*,S_{2h}}$.
For a failure probability $u\in(0,1)$, define the relative
observation-probability error scale
$\delta_m(u)=C\{\sqrt{\log(4JL/u)/(mh)}+\log(4JL/u)/(mh)\}$.

\textbf{Lemma 1 (core/outer sandwich).} If
\(mh\geq C\log(4JL/u)\) and \(\delta_m(u)\leq1/2\), then with probability
at least \(1-u\), the inclusions \(S_{2h}\subseteq\widehat S_{h,f}\subseteq S_{h/2}\) hold simultaneously over all folds, and every item in the outer set has relative observation-probability error at most \(\delta_m(u)\). See Appendix~\ref{app:spectral} for the proof.

\paragraph{Initialization under data-dependent preprocessing.}
\label{sec:residual}

The dependence created by data-dependent preprocessing is handled by an
exact pathwise decomposition. For a training set $T$ of size $m$, condition
on its indices and fixed labels, and set
$c_j=m^{-1}\sum_{i\in T}\theta_{z_i j}$. On the deterministic outer set
$S_o=S_{h/2}$, define the oracle residual and centered-mask matrices
$\mat R_\gamma$ and $\mat M_\gamma$, together with a diagonal
selection/rescaling matrix $\mat D_\gamma$ whose off-diagonal entries are zero:
\begin{align*}
R_{\gamma,ij}&=\pi_j^{-\gamma}\Omega_{ij}(X_{ij}-c_j),\\
M_{\gamma,ij}&=\pi_j^{-\gamma}(\Omega_{ij}-\pi_j),\\
D_{\gamma,jj}&=\mathbf1\{\widehat\pi_j\ge h\}
 \left\{\pi_j/\max(\widehat\pi_j,h)\right\}^{\gamma}.
\end{align*}

With discarded columns embedded as zeros, let $\mat H$ denote training-row
centering $\mat H_T$, suppress training subscripts, and set $\mat P_\gamma=\mat H(E\mat R_\gamma)\mat D_\gamma$.
The following identity holds on this outer set:
\begin{align}
\widehat{\mat R}_\gamma
 &=\mat P_\gamma+\mat H(\mat R_\gamma-E\mat R_\gamma)\mat D_\gamma
 \notag\\
 &\quad+\mat H\mat M_\gamma\mat D_\gamma
       \operatorname{diag}(\boldsymbol c-\widehat{\boldsymbol c}).
\label{eq:preprocessing-decomposition}
\end{align}

Lemma~1 places the actually retained columns inside this outer set with high
probability. The reference signal $\mat P_\gamma$ remains random through
$\mat D_\gamma$; the identity does not replace it with the expectation of the
empirically transformed matrix. Concentration is applied to the oracle noise
and mask matrices before multiplication by $\mat D_\gamma$ and the center
error. Simultaneous bounds and pathwise norm inequalities then bound the
resulting products without assuming independence between the estimates.
At $\gamma=0$, $\mat D_0$ remains a random selection projection.

\paragraph{Perturbation and separation scales.}
Fix $\gamma\in[0,1/2]$, $h\in(0,1/2]$ and $u\in(0,1/2)$.
For $m$ training rows and $J$ items, put $S_o=S_{h/2}$,
$S_c=S_{2h}$ and the log factor $\ell_u=\log\{C(m+J)K/u\}$.
The relative observation-probability error allowance is
$\delta=C\{\sqrt{\ell_u/(mh)}+\ell_u/(mh)\}$.
Here $v_\gamma$ measures the allocation of noise across items;
$B_\gamma$ and $B_\gamma^{\rm mask}$ bound response and mask noise, while
$a_c$ bounds the item-center estimation error. These quantities combine into
the perturbation scale $b_\gamma$; the corresponding core center separation
is $\Delta_{\gamma,c}$:
\begin{align*}
v_\gamma&=\max\left\{m\max_{j\in S_o}\pi_j^{1-2\gamma},
                         \sum_{j\in S_o}\pi_j^{1-2\gamma}\right\},\\
B_\gamma&=C\{\sqrt{\rho v_\gamma\ell_u}+h^{-\gamma}\ell_u\},\qquad
B_\gamma^{\rm mask}=C\{\sqrt{v_\gamma\ell_u}+h^{-\gamma}\ell_u\},\\
a_c&=C\left\{\sqrt{\frac{\rho\ell_u}{mh}}+\frac{\ell_u}{mh}\right\},\qquad
b_\gamma=(1-\delta)^{-\gamma}(B_\gamma+a_c B_\gamma^{\rm mask}),\\
\Delta_{\gamma,c}^2&=\min_{k\ne\ell}\sum_{j\in S_c}
                      \pi_j^{2-2\gamma}(\theta_{kj}-\theta_{\ell j})^2.
\end{align*}

\Needspace{9\baselineskip}
\textbf{Theorem 1 (initialization under data-dependent preprocessing).}
\label{prop:initialization}
Consider $m$ training rows and $J$ items under A1--A3, with training
class fractions bounded below. Set the rounding factor to $a_K=8(\log K+2)$.
Use the scales just defined and suppose
$\delta\le1/2$ and $\Delta_{\gamma,c}>0$.
Use empirical observation probabilities, Jeffreys item centers, the leading
$(K-1)$ scores $\mat{U}\mat{\Sigma}$ (with zero padding when needed), the
return and extended-loss conventions above, and $R$ independent D-squared
starts followed by nonincreasing Lloyd updates. For every tail multiplier
$\lambda>1$, with probability at least
$1-Cu-\lambda^{-R}$,
\[
\operatorname{Mis}_\gamma\le
\min\left\{1,\frac{C a_K\lambda(K-1)b_\gamma^2}{m\Delta_{\gamma,c}^2}\right\}.
\]

Taking $\lambda=2$ yields the fixed-factor bound with
$A_K=16(\log K+2)$ and exceptional probability $Cu+2^{-R}$.
Constants depend on balance and regular-profile bounds. The proof is in Appendix~\ref{app:residual-recovery}.

The class-center separation condition on the core items accommodates
centered signals of rank below $K-1$. The weighting exponent $\gamma$
affects both the preprocessing perturbation bound $b_\gamma$ and the
class-center separation $\Delta_{\gamma,c}$ on the core items. With other
conditions fixed, a smaller ratio $b_\gamma/\Delta_{\gamma,c}$ gives a
smaller error term in the training misclassification bound.
 
\paragraph{From training labels to likelihood recovery.}
The initialization condition connects training loss to held-out classification.
For a deterministic loss bound $\epsilon_m\in[0,1]$ and failure
probability bound $\eta\in(0,1)$, assume that
\begin{equation*}
\Pr\left\{\bigcap_{f=1}^L
\left[\operatorname{err}(\widehat{\boldsymbol z}^{(f)},\boldsymbol z_{T_f})
\le\epsilon_m\right]\right\}\ge1-\eta.
\tag{Init}\label{eq:initialization}
\end{equation*}

Each initializer is measurable with respect to the sigma-field generated by
its training information:
\[
\mathcal T_f=\sigma\{\text{split},
 (\Omega_{ij},\Omega_{ij}X_{ij})_{i\in T_f,j\le J},
 \text{training seeding for fold }f\}.
\]

In the proof, after successful initialization in fold $f$, we take the first
permutation $\sigma_f$ in a fixed order that minimizes training misclassification.
With true labels fixed, $\sigma_f$ is determined by that fold's training
information, so we can condition on this information when bounding the
misclassification probability for a held-out row. The algorithm aligns class
numbers across folds by matching their estimated profiles. Under the conditions
below, these class numbers correspond to the true classes through a single
global permutation.

In $r_m$, the first term accounts for contamination from the training
partition; the remaining terms describe sampling and smoothing error. The
prior-error scale is $a_m$, while $D_{h/2}$ and $y_n$ describe the observed
row counts that convert profile error into score error. Define these scales;
the superscript $J$ on $R_m^J$ denotes Jeffreys smoothing:
\begin{align*}
r_m&=\frac{\epsilon_m}{h\rho}
+\sqrt{\frac{\log(8KJL/\eta)}{mh\rho}}
+\frac{\log(8KJL/\eta)}{mh\rho},\\
a_m&=\epsilon_m+m^{-1},\qquad
D_{h/2}=\rho\sum_{j\in S_{h/2}}\pi_j,\\
y_n&=\log(8nKL/\eta),\\
R_m^J&=C\left[a_m+r_m\{D_{h/2}+\sqrt{D_{h/2}y_n}\right.\\*
&\hspace{10em}\left.{}+y_n\}\right].
\end{align*}

Profile matching compares the squared Hellinger separation
$D^{\mathrm{prof}}_{*,2h}$ with the estimation budget $e_m$.
Using the distance defined in Section~\ref{sec:method}, set
\[
D^{\mathrm{prof}}_{*,2h}=\min_{k\ne\ell}d_{S_{2h}}^2(\boldsymbol{\theta}_k,\boldsymbol{\theta}_\ell),
\qquad e_m=C|S_{h/2}|\rho r_m^2.
\]

On the simultaneous profile-estimation and item-screening events, the
foldwise Hungarian match is uniquely population-correct under the matching
condition
\begin{equation*}
16e_m<D^{\mathrm{prof}}_{*,2h}.
\tag{M}\label{eq:matching}
\end{equation*}

For a score-error budget $\xi\ge0$, let $u=1-s\in[0,1]$ be the
optimizing variable and define the penalized Chernoff exponent on $S$ by
\begin{equation*}
\begin{aligned}
\mathcal C_{k\ell,S}(\xi)
&=\sup_{0\leq u\leq1}\left\{
-\sum_{j\in S}\log c_{k\ell j}(1-u)-2u\xi
\right\},\\
\mathcal C_{*,S}(\xi)
&=\min_{k\ne\ell}\mathcal C_{k\ell,S}(\xi).
\end{aligned}
\end{equation*}

\Needspace{9\baselineskip}
\textbf{Theorem 2 (cross-fitted likelihood refinement).}
Assume the model and split conditions A1--A4, the initialization
condition~\eqref{eq:initialization}, and the profile-matching
condition~\eqref{eq:matching}. For sufficiently small fixed constants
$c_0,c>0$, suppose
\begin{equation*}
\begin{aligned}
mh&\geq C\log(4JL/\eta),&
mh\rho&\geq C\log(8KJL/\eta),\\
\delta_m(\eta)&\leq c_0,\quad Cc_0\leq\tfrac12,&
r_m&\leq c.
\end{aligned}
\end{equation*}

Then the Jeffreys-profile, Hellinger-matched, cross-fitted observed-likelihood
reassignment satisfies, up to one global permutation,
\begin{equation}
P_{\mathrm{perm}}\{\widehat{\boldsymbol{z}}\ne \boldsymbol{z}\}
\leq C\eta+C n(K-1)\exp\{-\mathcal C_{*,S_{2h}}(R_m^J)\}.
\label{eq:recovery-bound}
\end{equation}

At $\xi=0$, \(\mathcal C_{k\ell,S}(0)\) equals the oracle observed Chernoff
information. More generally,
\begin{equation*}
\mathcal C_{k\ell,S}(\xi)
\geq C^{\mathrm{obs}}_{k\ell,S}-2\xi,
\qquad
\mathcal C_{k\ell,S}(\xi)
\geq \tfrac12I^{\mathrm{obs}}_{k\ell,S}-\xi.
\end{equation*}

The second inequality evaluates the penalized exponent at the feasible value
$s=1/2$, and therefore uses one half of the order-$1/2$ R\'enyi divergence.
Appendix~\ref{app:refinement} proves the recovery bound.
Under the model, matching condition, and concentration conditions of
Theorem~2, exact recovery follows if, for some fixed $\delta>0$,
\begin{equation*}
R_m^J=o(I^{\mathrm{obs}}_{*,2h}),\qquad
I^{\mathrm{obs}}_{*,2h}\geq(2+\delta)\log n,\qquad
\epsilon_m=o(h\rho),\qquad \eta=o(1).
\end{equation*}

\subsection{The sharp threshold for exact recovery}
\label{sec:sharp-recovery}
The hardest class pair determines the oracle Chernoff exponent for classifying
one subject. The union bound over subjects introduces a factor of $n$, so the
leading likelihood term is governed by $\log n-C_*^{\rm obs}$. To attain this
boundary with estimated profiles, we need the loss in the exponent to be
small relative to $C_*^{\rm obs}$. A matching converse must also account for
the possibility of assigning all labels jointly using the known class sizes.

Near this boundary, estimates of individual response probabilities can be
too noisy for the uniform bound in Theorem~2. Theorem~3 analyzes the same
estimator by bounding the aggregate effect of profile errors across items
on the Chernoff exponent.

\textbf{Theorem 3 (sharp first-order exact-recovery threshold).}
Let $K,L\ge2$ and the restart count $R=R_0\ge1$ be fixed, let $J\asymp n$,
and suppose that class proportions are bounded below. Fix an observation-probability
comparison factor $C_\pi$ and
a separation lower bound $d_0>0$. Assume A2--A4, $\pi\le\pi_j\le C_\pi\pi$,
$\theta_{kj}=\rho q_{kj}$, $0<q_-\le q_{kj}\le q_+$,
$\rho q_+\le1-\epsilon_0$, and
\[
\min_{k\ne\ell}J^{-1}\sum_j(q_{kj}-q_{\ell j})^2\ge d_0>0.
\]

Put $d_n=n\pi\rho$ and $C_*^{\rm obs}=\min_{k\ne\ell}C_{k\ell}^{\rm obs}$.
Use Algorithm~\ref{alg:crossfit} with $h=\pi/2$, any fixed
$\gamma\in[0,1/2]$, exact truncated $\mat U\mat\Sigma$ scores and $R_0$
independent standard D-squared starts followed by nonincreasing Lloyd updates.
For every fixed $a>0$, if $d_n\ge a\log n$, then for some fixed
$\beta\in(0,1/2]$ and $C<\infty$, for all sufficiently large $n$,
\begin{align}
P\{\operatorname{err}>0\}
&\le Cn^{-2}+CLd_n^{-R_0/2}\notag\\
&\quad+Cn(K-1)\exp\{-C_*^{\rm obs}
                 +C(\sqrt{d_n}+d_n^{1-\beta})\}.
\label{eq:sharp-recovery-bound}
\end{align}

Here $\beta$ keeps the oracle optimizing order away from zero and one.
Conversely, under the same model assumptions but without requiring
$d_n\ge a\log n$, fix a parameter sequence and balanced integer class
sizes, and let $\mathcal Z_n$ contain all label vectors of those sizes.
Write $P_{\boldsymbol z}$ for probability under a specified label vector.
Even when the estimator knows the class profiles, observation probabilities,
and class sizes, for any fixed $\delta>0$ we have
\[
C_*^{\rm obs}\le(1-\delta)\log n
\quad\Longrightarrow\quad
\sup_{\widehat{\boldsymbol z}}\inf_{\boldsymbol z\in\mathcal Z_n}
 P_{\boldsymbol z}\{\operatorname{err}=0\}\longrightarrow0.
\]

\paragraph{Proof of Theorem~3.}
For the upper bound, the training events in Appendix~\ref{app:sharp-outline}
give training loss $O(d_n^{-1/2})$
and at most $C/\pi$ discarded items per fold, with exceptional probability
at most $Cn^{-2}+CLd_n^{-R_0/2}$. The class profiles are estimated from
the same rows used to form the initial groups. The bound
$\|\boldsymbol w^\top\mat E\|_2^2\le \|\boldsymbol w\|_2^2\|\mat E\|^2$
holds for every group indicator $\boldsymbol w$, where $\mat E$ is the
centered training success matrix. Equation~\eqref{eq:sharp-aggregate}
in Appendix~\ref{app:sharp-outline} uses this bound to give squared
Hellinger error at most $CJ\rho/d_n$ for each fitted profile. Two folds'
retained sets omit at most $C/\pi=o(J)$ items from their intersection.
The true profiles have squared Hellinger separation at least $cJ\rho$
on that intersection. The triangle inequality then makes each correct
match a strict row minimum, so all folds share one label permutation.

For a held-out row $i$ of class $k$, let $\widehat L_{ik}^{(f)}$ be the
fitted log score for class $k$, using its fold's profiles relabeled by the
training permutation. Conditional on the corresponding training
information $\mathcal T_f$, Markov's inequality at the full-product oracle
order $u_*$ gives, on the training event for fold $f$ and for $\ell\ne k$,
\[
P_k\{\widehat L_{i\ell}^{(f)}\ge\widehat L_{ik}^{(f)}\mid\mathcal T_f\}
\le
(\widehat\alpha_{\ell}^{(f)}/\widehat\alpha_{k}^{(f)})^{u_*}
\prod_{j\in\widehat S_{h,f}}\widehat c_j(u_*).
\]

The fitted class proportions are bounded below. Each omitted item
contributes at most $C\pi\rho$ to the oracle exponent at $u_*$, so
screening costs at most $C\rho$. The fractional-moment bound in
Appendix~\ref{app:sharp-outline} therefore bounds the display by
$C\exp\{-C_{k\ell}^{\rm obs}+C(\sqrt{d_n}+d_n^{1-\beta})\}$.
Summing over rows and competitors and adding the probability that the
training events fail gives Equation~\eqref{eq:sharp-recovery-bound}.
By Proposition~3, $C_*^{\rm obs}\asymp d_n$. The exponent correction is
lower order, so the bound vanishes at a fixed positive margin above
$\log n$.

For the lower bound, first let $d_n\to\infty$. On either side of a hardest
class pair, the strict likelihood reversal bound in
Appendix~\ref{app:sharp-outline} gives probability at least
$c\exp\{-C_*^{\rm obs}-C\sqrt{d_n}\}$ that a row favors the other class.
Below $(1-\delta)\log n$, we have $d_n=O(\log n)$. Each class has order
$n$ independent rows, so the expected number of reversals in each
class is at least $n^{\delta-o(1)}$, and both classes contain a reversed
row with probability tending to one. Swapping those two labels preserves
the class sizes and increases the likelihood. Under the uniform prior on
$\mathcal Z_n$, the named-label MAP estimator therefore has success
probability tending to zero. Each permutation orbit contains at most $K!$
assignments, so the optimal Bayes success up to permutation also tends
to zero and bounds the minimax success from above. If $d_n$ stays bounded,
the disjoint-pair bound in Appendix~\ref{app:sharp-outline} gives the same
conclusion. Every
subsequence admits a further subsequence in one of these two regimes.

\textbf{Corollary 1 (common sharp threshold for residual and mean-imputation recovery).}
Under the model and algorithm conditions of Theorem~3, for any fixed $\delta>0$,
\begin{equation}
C_*^{\rm obs}\ge(1+\delta)\log n
\quad\Longrightarrow\quad P\{\operatorname{err}=0\}\longrightarrow1.
\label{eq:common-order-condition}
\end{equation}

This holds for every fixed $\gamma\in[0,1/2]$, including Mean and Residual,
and also for the column-centered, empirically observation-probability-weighted
zero-filled input
\begin{equation}
\widehat{\mat R}_{\mathrm{raw},T}
=\mat H_T(\mat\Omega\odot\mat X)_{S(T)}
 \operatorname{diag}(\widehat\pi_{S(T)}^{-1/2}).
\label{eq:centered-raw}
\end{equation}

In each case, use the prescribed effective dimension
$r_f=\min\{K-1,|\widehat S_f|,|T_f|-1\}$, with zero padding when needed.
The centered signal may have rank below $K-1$. Together with the converse,
this establishes the same first-order threshold for these representations
under fixed positive margins. The proof is in Appendix~\ref{app:common-threshold-proof}.

The common threshold applies to different initialization representations
followed by the same final reassignment. In every version, reassignment uses
the observed-data likelihood based on $\mat W$, which retains the distinction
between observed zeros and missing responses. Proposition~2 quantifies the
information lost when this distinction is removed. Section~\ref{sec:threshold-experiment}
examines how much recovery rates can vary at finite sample sizes.

For example, take $\pi_j=1$ and regular separated profiles with
$\rho\asymp\log n/n$ calibrated to $C_*^{\rm obs}/\log n=1+\delta$.
Corollary~1 applies, whereas the sufficient scaling from the uniform
Theorem~1--Theorem~2 route, $n\pi^2\rho^2/\log n\to\infty$, does not hold.
Appendix~\ref{app:route-comparison} records this substitution explicitly.
The threshold $h=\pi/2$ is calibrated to the specified comparable
observation-rate scale.

\section{Informative-core extensions}
\label{sec:extensions}
A deterministic informative core can support exact recovery when other items
are much rarer or lie at the screening boundary. Corollary~2 combines
exclusion of the rare tail with the core-only sharp threshold of Corollary~1.
Corollary~3 applies the general bounds of Theorems~1--2 at a stronger
uniform-profile scale and allows the retained item set to remain random.

\subsection{Recovery on an informative core}
\label{sec:stable-core}

\textbf{Corollary 2 (informative core with an arbitrarily rare tail).}
Partition the items into a deterministic core $S_c$ of size
\(J_c=|S_c|\asymp n\) and a set of rarely observed items, called the tail,
of size $T_n=J-J_c$.
Assume A1--A2 and A4, with fixed $K,L$. Let $p_n>0$ be the core
observation-probability scale. On the core, assume the regularity and
separation conditions of Theorem~3, with \(p_n\le\pi_j\le C_\pi p_n\), and set
\(h=p_n/2\). On the tail, allow arbitrary Bernoulli profiles and
\(0\le\pi_j\le p_n/8\). Write
$C_{*,S_c}^{\rm obs}=\min_{k\ne\ell}C_{k\ell,S_c}^{\rm obs}$ for the
minimum core Chernoff information. Suppose that for some fixed $\delta>0$,
\[
C_{*,S_c}^{\rm obs}\ge(1+\delta)\log n,
\qquad \log(T_n+1)=o(n p_n).
\]

Under these conditions, the cross-fitted estimator using either mean-imputation or
observation-probability-weighted residual scores, with the algorithm
conventions of Theorem~3 and any fixed $R_0\ge1$ restarts, recovers all
labels up to one permutation with
probability tending to one. The proof is in Appendix~\ref{app:rare-tail-proof}.

\subsection{Recovery without retained-set stabilization}
A sufficiently informative core supports recovery when items at the
screening boundary are retained at random.

\textbf{Corollary 3 (recovery with persistently random item screening).}
\label{cor:random-retention}
Let $K,L$ and $R_0\ge1$ be fixed, $J\asymp n$, and assume A1--A4.
For a deterministic core $S_c$ with $|S_c|\asymp n$, suppose that for
$0<h\le1/2$ and a fixed $C_\pi$,
\[
2h\le\pi_j\le C_\pi h\quad(j\in S_c),\qquad
0\le\pi_j\le C_\pi h\quad(j\notin S_c).
\]

Under the regular profile bounds,
assume
\[
\min_{k\ne\ell}|S_c|^{-1}\sum_{j\in S_c}(q_{kj}-q_{\ell j})^2\ge d_0>0,
\qquad \frac{nh^2\rho^2}{\log n}\longrightarrow\infty.
\]

Then the same cross-fitted estimator, for every fixed
$\gamma\in[0,1/2]$, recovers all labels up to one permutation with
probability tending to one. The proof is in Appendix~\ref{app:random-retention-proof}.

No gap around $h$ is required for noncore observation probabilities. In
particular, a positive fraction of items may satisfy $\pi_j=h$. Because
$mh\to\infty$ and $h\le1/2$, the binomial central limit theorem implies an
asymptotic screening probability of $1/2$ for each such item in a training
fold. Thus recovery of the full retained item set is unnecessary for label
recovery. The core belongs to every retained intersection with high
probability, while the randomly screened boundary items contribute to the
noise and score budgets through the outer set.

\section{Simulation}
\label{sec:simulation}

We first compare finite-sample recovery across observed-information levels
and profile sources near the threshold in Theorem~3. We then examine screening
and rounding, and assess spectral initialization for iterative likelihood
fitting. We report the fraction of datasets in which a method recovers all
labels up to a global permutation. Paired comparisons use independent datasets as the units for
Monte Carlo standard errors (MCSEs) \citep{MorrisEtAl2019}. Failed returns
receive loss one, and capped EM fits enter the summaries with their returned
outcomes.

\subsection{Recovery near the observed-information threshold}
\label{sec:threshold-experiment}

Both designs have two equally sized classes and calibrate
$c=C_*^{\rm obs}/\log n$. The symmetric design alternates class profiles
across items and sets $c=.7$ or $1.3$ at $n=J=240,480,960$. Three paths
vary response sparsity, observation sparsity, or both. The asymmetric design
uses three item types with different observation probabilities and asymmetric
class profiles. It sets $c\in\{.7,.9,1,1.1,1.3\}$ at $n=J=480,960$.
Each cell uses 100 datasets. The fitted methods use two folds, screen at
half the minimum observation probability, and perform one Jeffreys likelihood
reassignment. The rowwise oracle knows the profiles and equal class priors.
The count oracle also enforces the true class sizes.

\begin{figure}[!htbp]\centering
\includegraphics[width=\textwidth]{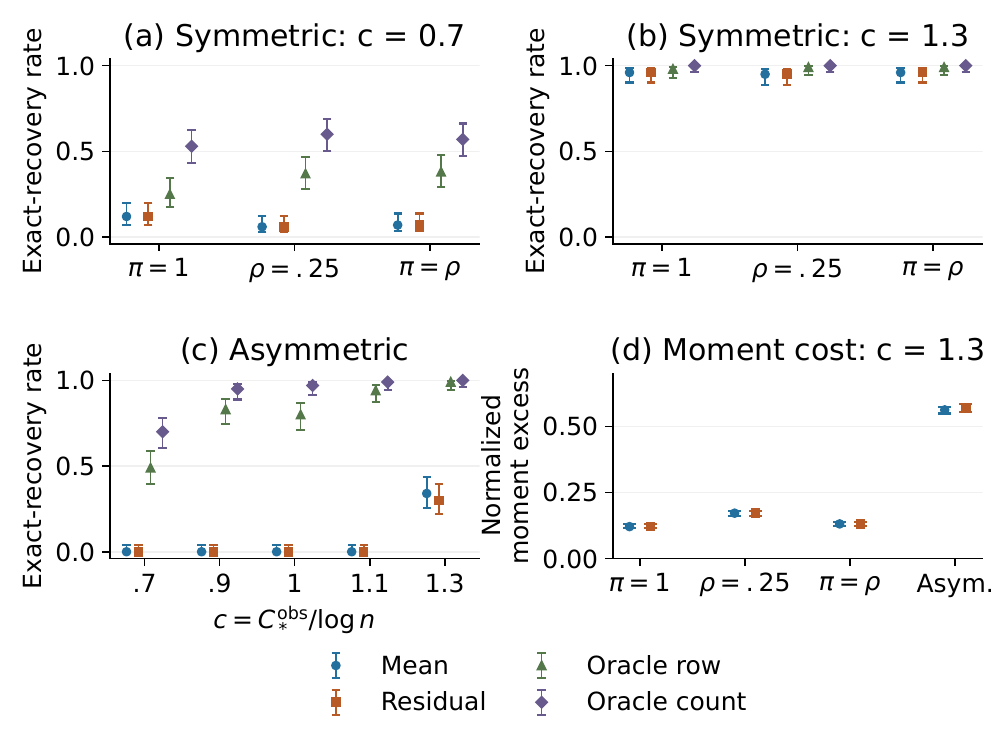}
\caption{Observed information and finite-sample recovery at $n=J=960$.
(a)--(c) show exact-recovery rates with pointwise 95\% Wilson intervals from
100 datasets per cell. Each fitted method uses ten D-squared starts per fold
and one cross-fitted likelihood reassignment. Panel (d) shows the median and
interquartile range of $M$. For each dataset, $M$ is the maximum over two
folds and ordered class pairs of the conditional plug-in log-moment excess
on retained items at the oracle order, divided by full oracle $C_*^{\rm obs}$.
The diagnostic excludes priors.}
\label{fig:threshold-comparison}
\end{figure}

Figure~\ref{fig:threshold-comparison} shows how much recovery can differ at
the same value of $c$. Across the symmetric paths and sizes, the ten-start
fits achieve exact recovery in 88--98 of 100 datasets at $c=1.3$ and in
4--24 at $c=.7$. At $n=960$ and $c=1.3$, the symmetric complete-observation
path has a median maximum coordinatewise profile error of $3.57\rho$ and
a median $M$ of $.121$ (Online Resource~1, Table~\ref*{tab:sharp-diagnostics}).
In the asymmetric design, the fitted methods recover all labels in 21--34 of
100 datasets at $c=1.3$. At both sizes, the rowwise oracle achieves exact
recovery in 99 of 100 datasets and the count oracle in all 100. At the tested
values $c\le1.1$, neither fitted method recovers all labels in any dataset.
At $n=960$ and $c=1.3$, with class priors estimated within each training fold,
Mean recovers all labels in 34 datasets and Residual in 30. They misclassify
an average of 1.20 and 1.31 of the 960 subjects, respectively.

On the same 100 asymmetric datasets at $n=J=960$ and $c=1.3$, we compare
four profile sources using the same folds and retained items, with class
priors fixed at $1/2$. The rowwise classifier recovers all labels in 99 datasets
with the true profiles and in 52 with profiles estimated from the true training labels, a
paired drop of 47 percentage points (MCSE 5.02 percentage points). Using
the Mean and Residual spectral groups to estimate the profiles yields 36
and 32 exact recoveries, respectively. The median $M$ across datasets is
$0$, $.493$, $.562$, and $.570$ in the same order (Online Resource~1,
Section~\ref*{sec:asymmetric-threshold-diagnostics}, Table~\ref*{tab:paired-profile-cost}).
The largest drop in exact recovery occurs when the profiles are estimated
within the correct training groups.

\subsection{Screening and spectral rounding}
\label{sec:recovery-mechanism}
\label{sec:narrative-recovery}

We examine the informative-core mechanism in Corollary~3 using three classes
with vanishing response and observation probabilities and a rank-one centered
signal. A quarter of the items sit at the screening threshold. Across
$n=J=240,480,960$, every core item is retained in the recorded folds, while
roughly half the boundary items are selected. With ten starts, Mean and
Residual each recover all labels in 94 of 100 datasets at the smallest size
and in all 100 at each larger size. With one start at the largest size,
Mean recovers all labels in 91 of 100 datasets and Residual in 82.

\begin{figure}[!htbp]\centering
\includegraphics[width=\textwidth]{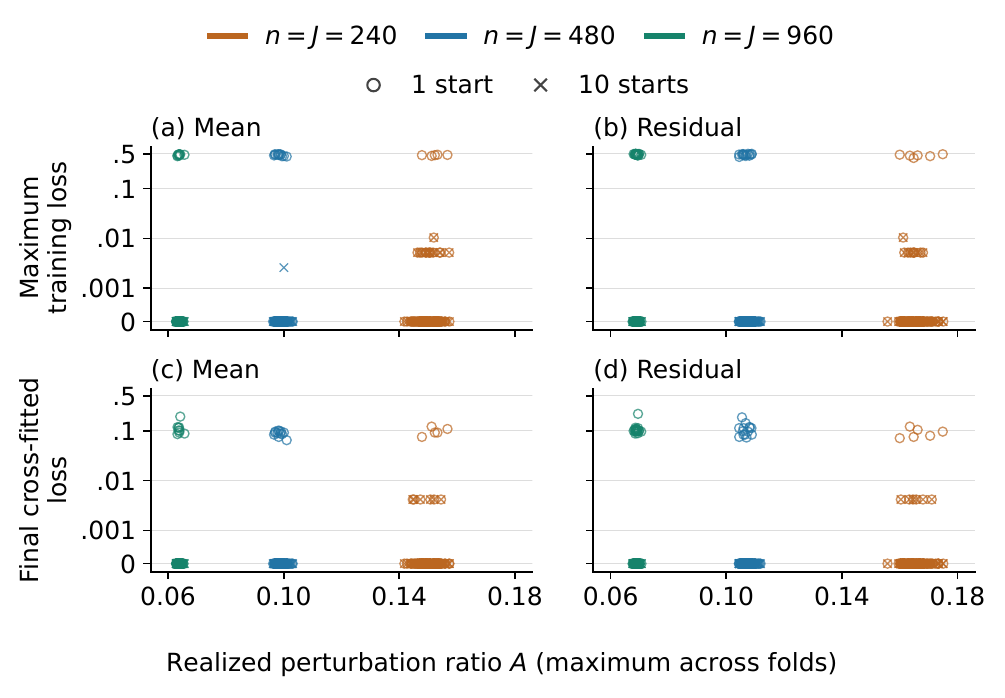}
\caption{Realized perturbation and classification loss on the sparse
random-screening path. The horizontal coordinate is the maximum over five
folds of $A=(K-1)\|\widehat{\mat R}_\gamma-\mat P_\gamma\|^2/(m\Delta_{\gamma,c}^2)$.
Here $\|\cdot\|$ is the operator norm.
The top row shows maximum training classification loss and the bottom row
shows final cross-fitted loss. One and ten starts share each preprocessed
matrix. Each size contributes 100 datasets per start budget, and points can
overlap. The vertical scale is linear near zero and logarithmic above .001.}
\label{fig:perturbation-recovery}
\end{figure}

In this design, the median perturbation ratio falls as the numbers of
subjects and items both increase from 240 to 960. For a given dataset and
spectral representation, the one- and ten-start runs use the same spectral
scores, so the differences in classification loss between these runs in
Figure~\ref{fig:perturbation-recovery} arise from rounding. Theorem~1 uses
the perturbation ratio to bound training misclassification, and the number
of starts enters the upper bound on the probability of rounding failure.
\FloatBarrier

\subsection{Spectral initialization for likelihood fitting}
\label{sec:em-budget}

This comparison uses iterative observed-data EM and evaluates classification
loss on the full training sample. Four simulation families (A8, B4, D3, and
M5) each contribute 40 paired datasets of 800 subjects. Spectral initializers
use ten D-squared starts and one EM trajectory. Random~1 and Random~10 use
one and ten EM trajectories from random parameter perturbations. Random~10
selects the returned fit with the largest training likelihood. Each trajectory
uses a cap of 150 iterations and a relative log-likelihood tolerance of
$10^{-6}$, with the updates given in Appendix~\ref{app:em-details}.

In the five-class M5 family, the spectral starts yield lower mean classification
loss and shorter median fitting times than both random-start procedures
(Table~\ref{tab:narrative-em-budget}). Final losses are similar across all four
procedures in the other three families.

\par\medskip\noindent\begin{minipage}{\textwidth}\centering
\makeatletter\def\@captype{table}\makeatother
\caption{Classification loss and computing time for EM on 40 paired M5 datasets. Loss is mean (MCSE), time is median [IQR] in seconds and includes initialization and all attempted EM trajectories, and NC counts selected fits that did not meet the convergence criterion. All returns enter the summaries.}
\label{tab:narrative-em-budget}
\begin{tabular}{@{}lrrrr@{}}
\toprule
Start & EM trajectories & Loss (MCSE) & Seconds [IQR] & NC\\
\midrule
Mean & 1 & 0.1557 (0.0020) & 0.53 [0.46,0.62] & 0\\
Residual & 1 & 0.1588 (0.0021) & 0.51 [0.39,0.66] & 0\\
Random 1 & 1 & 0.2350 (0.0156) & 1.44 [1.25,1.83] & 3\\
Random 10 & 10 & 0.1680 (0.0026) & 17.19 [16.07,18.29] & 0\\
\bottomrule\end{tabular}
\end{minipage}\par\medskip

In the paired M5 comparison, we hold the screening and EM update rule fixed
and vary the spectral representation. Switching from Mean to Residual reduces
initial loss by .0082 (MCSE .0033) and increases final loss by .0031
(MCSE .0011, Online Resource~1, Section~\ref*{sec:preprocessing-b},
Table~\ref*{tab:preprocessing-b-paired}).
The ranking of the two spectral starts by classification loss reverses after EM.
\FloatBarrier

\section{Real-data analysis}
\label{sec:pisa}
In the PISA 2022 mathematics data, we compare how models fitted from
different EM starts predict the observed responses of students in held-out
schools. After screening, we retained 154 common dichotomous items and 4,309
students in 154 United States (US) schools and 6,167 students in 164
Singapore schools.
Students have scores on an average of 18.91 items in the USA and 18.93 in
Singapore. Adaptive item routing \citep[Chapter~2]{OECD2024} and nonresponse place the observation process outside
the independent-MCAR model. The students' class labels are unknown.

We grouped US schools into five folds and used the one-standard-error (one-SE)
rule to select $K=4$ and $\tau=0$. Under this rule, we choose the simplest
candidate whose score is within one standard error of the best score. The prespecified
validation-gated selector uses the propensity-weighted Raw representation
with centered rank-$(K-1)$ $\mat U$ scores. The EM comparisons share the
selected $K$ and $\tau$. Mean and Residual each use their own representation
with $\mat U\mat\Sigma$ scores. For Singapore, we kept the US settings and
the same 154 scoring items and re-estimated the model parameters using the
training schools in each fold. The US scores are post-selection. We used
data from both countries to screen the common item pool.

For each held-out student, the predictive score is the whole-row marginal
log likelihood divided by the number of scored items. Students receive equal
weight. Both comparators use the unit-slope Rasch model \citep{Rasch1960}.
We assume a zero-mean normal ability distribution, fixing its variance at one
for Rasch and estimating it for Rasch-SD.

MeanEM and ResidualEM use one EM trajectory after ten spectral-rounding
starts. The random comparators use one or ten EM starts. At the common
150-iteration cap, the Residual-minus-Mean score is $+.00032$ (standard error .00038)
in the USA and $+.00015$ (.00012) in Singapore (Table~\ref{tab:focused-pisa}).
Continuing the seven capped single-trajectory fits to convergence changes
the US contrast to $-.00011$ (.00029). The Singapore contrast stays at
$+.00015$ (.00012), where both spectral fits had already converged in
every fold (Online Resource~1, Section~\ref*{app:pisa-diagnostics}).

\par\medskip\noindent\begin{minipage}{\textwidth}\centering
\makeatletter\def\@captype{table}\makeatother
\caption{Paired PISA marginal-score differences at the 150-iteration EM cap, using the same held-out students and a common pool of 154 items. Larger values are better. Parentheses give descriptive standard errors clustered by school.}
\label{tab:focused-pisa}
\begin{tabular}{@{}lrr@{}}
\toprule
Contrast & USA & Singapore\\
\midrule
Residual $-$ Mean & +0.0003 (0.0004) & +0.0001 (0.0001)\\
Residual $-$ Random 1 & -0.0009 (0.0003) & +0.0007 (0.0002)\\
Residual $-$ Random 10 & -0.0010 (0.0003) & 0.0000 (0.0001)\\
Mean $-$ Rasch-SD & +0.0076 (0.0009) & +0.0038 (0.0005)\\
Residual $-$ Rasch-SD & +0.0079 (0.0009) & +0.0040 (0.0005)\\
Rasch-SD $-$ Rasch & +0.0001 (0.0001) & +0.0020 (0.0003)\\
\bottomrule
\end{tabular}
\end{minipage}\par\medskip

\section{Discussion}
\label{sec:discussion}

We found that a cross-fitted classifier can attain the oracle's first-order
threshold for exact label recovery even when it learns both its spectral
representation and its class profiles from incomplete responses. In each
training fold, we used the same subjects to form the groups and estimate
their response probabilities. Theorem~3 shows that the resulting profile
errors, taken across items, reduce the classification error exponent by only
a lower-order term. Both Mean and Residual then refine the labels using a
likelihood that distinguishes a missing response from an observed zero.
Our paired adaptation of SOLA \citep{LyuGu2026} illustrates why this
distinction matters. We kept the spectral start fixed and found that using the
observation mask during refinement reduced mean classification loss in every simulation
family (Online Resource~1, Table~\ref*{tab:focused-sola-paired}).

The lower-order term can still matter at the sample sizes we studied. In the
asymmetric experiment at $n=J=960$ and $C_*^{\rm obs}=1.3\log n$, the
rowwise classifier recovered every label in 99 of 100 datasets with the true
profiles. When we estimated the profiles from the true training groups,
keeping the folds, retained items, and equal class priors fixed, the count fell to 52.
Using spectral groups caused a further, smaller decline
(Section~\ref{sec:threshold-experiment}). Even a correct training group
contained few observed successes on some items, leaving their response
probabilities difficult to estimate. Thus the largest loss in this design
arose from estimating the profiles within groups. A common first-order
threshold allows a substantial finite-sample gap when a single incorrect
label suffices to prevent exact recovery.

Our proof bounds the order of the estimation cost. A second-order analysis
could quantify how much more information the fitted classifier needs than
the oracle to achieve a specified recovery probability at a given sample
size. To test that prediction, we could vary class separation at several sample
sizes and use paired simulations to estimate the gap in information needed
for the two classifiers to reach the same exact recovery rate.
PISA poses a further difficulty because we observe students' responses
but lack their true class labels. In both countries, the
latent class models we fitted from the two spectral starts scored higher than
Rasch-SD on students in held-out schools. Yet the small US difference in
mean log score per scored item (Residual minus Mean) changed from $+0.00032$
at the iteration cap to $-0.00011$ at convergence (Section~\ref{sec:pisa}).
The predictive score averages response-pattern probabilities over classes,
so two fits can predict equally well while assigning different labels.
Item routing also makes the observation indicators dependent. Extending
the recovery result to this setting would require
deriving the information that distinguishes classes under the actual item
assignment design. We could then apply the PISA routing rules to simulated
students with known labels and measure how often the classifier recovers all
labels on either side of the predicted boundary.

Under the independent item-specific MCAR model and the remaining conditions
of Theorem~3, the minimum observed Chernoff information determines the
first-order boundary for recovering all $n$ labels. Above $\log n$ by a fixed
relative margin, Mean and Residual each recover every label up to permutation
with probability tending to one after a single cross-fitted likelihood update. Below it by a
fixed relative margin, even an oracle that knows the parameters and class
sizes has vanishing worst-case success probability. An imperfect training
partition can therefore suffice for exact recovery when its fitted profiles
preserve the evidence separating the hardest pair of classes.

\FloatBarrier
\appendix
\section{Information identities and implementation conventions}
\label{app:bridge}
The information scale used throughout the recovery argument is determined by
the observed three-outcome distribution. This appendix proves
Propositions~1--3 and identifies the information lost by zero filling.

\subsection{Observed information and zero filling}
\label{app:observed-information}
\paragraph{Proof of Proposition 1.}
Summing the three item atoms yields $a^{\rm obs}_{k\ell j}$ and
$c_{k\ell j}(s)$ for $0<s<1$, including boundary response probabilities.
Mutual independence of the cells multiplies these item coefficients.
Taking logarithms establishes Proposition~1, and taking the infimum over
$0<s<1$ yields the Chernoff formula in Section~\ref{sec:chernoff-information}.
When $b_{k\ell j}\to1$, expanding the logarithm gives the leading term
$2\pi_j(1-b_{k\ell j})$ for the observed item R\'enyi divergence, which is
$\pi_j$ times the complete-data leading term.

\paragraph{Proof of Proposition 2.}
The likelihood retains the full three-outcome symbol
\(\mat{W}=(\mat{\Omega},\mat{\Omega}\odot \mat{X})\), whereas the raw
spectral input retains only \(\mat{V}=\mat{\Omega}\odot \mat{X}\). For one item,
write $p=\theta_{kj}$, $q=\theta_{\ell j}$, $\pi=\pi_j$, and
$V=V_{ij}$. Denote the affinity of the zero-filled symbol by
$a^+(p,q,\pi)$. Its two class-conditional laws are
\begin{align*}
V\mid k&\sim\operatorname{Ber}(\pi p),\\
a^+(p,q,\pi)&=\pi\sqrt{pq}+\sqrt{(1-\pi p)(1-\pi q)},\\
I^+_{k\ell,S}&=-2\sum_{j\in S}
\log a^+(\theta_{kj},\theta_{\ell j},\pi_j).
\end{align*}

For the observed symbol, the missing atom is common to both classes. Thus,
canceling that atom and rationalizing the remaining difference gives
\begin{align*}
a^{\mathrm{obs}}(p,q,\pi)
&=1-\pi+\pi\{\sqrt{pq}+\sqrt{(1-p)(1-q)}\},\\
a^+-a^{\mathrm{obs}}
&=\frac{\pi(1-\pi)
\{\sqrt{1-p}-\sqrt{1-q}\}^{2}}
{\sqrt{(1-\pi p)(1-\pi q)}+(1-\pi)
+\pi\sqrt{(1-p)(1-q)}}\geq0.
\end{align*}

When \(\pi=1\) and \(\max\{p,q\}=1\), the displayed quotient is interpreted
by its continuous extension, which equals zero.

The numerator vanishes exactly when $p=q$, $\pi=0$, or $\pi=1$; continuity
handles the zero denominator at the boundary. Since $a^+\geq a^{\mathrm{obs}}$
and $-\log$ is decreasing, the R\'enyi ordering follows. The zero-filled
symbol is a deterministic coarsening of the observed symbol, so the same
ordering holds for Chernoff information by data processing. Under regular
scaling, uniformly for $p=\rho a$ and $q=\rho b$,
\[
(\sqrt{1-p}-\sqrt{1-q})^2
=\tfrac14\rho^2(a-b)^2+O(\rho^3(a-b)^2).
\]

The denominator in the rationalized expression is $2+O(\rho)$. Therefore,
\[
a^+-a^{\mathrm{obs}}
=\frac{\pi(1-\pi)\rho^2(a-b)^2}{8}
+O\!\left(\pi(1-\pi)\rho^3(a-b)^2\right).
\]

Both affinities equal $1-O(\pi\rho)$. Expanding $-2\log$ gives the
one-item expansions
\begin{align*}
I^{\mathrm{obs}}-I^+
&=\frac{\pi(1-\pi)}4\rho^2(a-b)^2
 +O\!\left(\pi(1-\pi)\rho^3(a-b)^2\right),\\
I^{\mathrm{obs}}
&=\pi\rho(\sqrt a-\sqrt b)^2+O\!\left(\pi\rho^2(a-b)^2\right).
\end{align*}

On the regularity interval,
$(a-b)^2\asymp(\sqrt a-\sqrt b)^2$. Summing the resulting nonnegative
itemwise bounds proves the relative inequality.

\subsection{Information equivalence under regular probability scaling}
\label{app:information}

\paragraph{Proof of Proposition 3.}
For one item, again set $p=\theta_{kj}=\rho a$ and
$q=\theta_{\ell j}=\rho b$. Its complete-response Bernoulli affinity is
$b(p,q)=\sqrt{pq}+\sqrt{(1-p)(1-q)}$, and hence
\[
1-b(p,q)=\frac12\left\{(\sqrt p-\sqrt q)^2
+(\sqrt{1-p}-\sqrt{1-q})^2\right\}.
\]

Regularity implies
$c_0\rho(a-b)^2\leq 1-b(p,q)\leq C_0\rho(a-b)^2$.
Write the observed item affinity as $1-x$, where
$x=\pi_j\{1-b(p,q)\}$. The upper probability bound gives
$0\leq x\leq 1-\epsilon_0$, and therefore
\[
2x\leq -2\log(1-x)\leq \frac{2x}{\epsilon_0}.
\]

Each observed R\'enyi contribution is consequently comparable to
$\pi_j\rho(a-b)^2=\pi_j(p-q)^2/\rho$. Summing over items proves the
first claim.

At the midpoint order $s=1/2$, the same comparison supplies a lower bound
of the order of $Q_{k\ell}$. For the upper bound, Chernoff information is
no larger than either directed Kullback--Leibler (KL) divergence,
$D_{\rm KL}(P\|Q)=\sum_xP(x)\log\{P(x)/Q(x)\}$, for discrete laws.
Let $P_k^{(j)}$ be the class-$k$ law of observed item $W_{ij}$.
The common missing atom contributes zero, and
\[
D_{\mathrm{KL}}(P_k^{(j)}\|P_\ell^{(j)})
=\pi_jD_{\mathrm{KL}}\{\operatorname{Ber}(p)\|\operatorname{Ber}(q)\}
\leq \pi_j\frac{(p-q)^2}{q(1-q)}
\leq C\pi_j\rho(a-b)^2.
\]

Applying this bound itemwise and using additivity of KL divergence for
product distributions yields the Chernoff upper bound. Together with the
midpoint lower bound, this proves $C^{\mathrm{obs}}_{k\ell}\asymp Q_{k\ell}$.

\subsection{Exact-arithmetic conventions}
\label{app:exact-conventions}
Algorithm~\ref{alg:crossfit} uses the leading $r_f$ scores, padding with
zeros to $K-1$ coordinates when needed. At a tied truncation boundary,
choose the retained subspace by a fixed rule based on the row Gram matrix,
such as projecting the row-coordinate basis onto the tied eigenspace and
orthonormalizing in index order. If D-squared sampling mass vanishes before
$K$ centers are selected, fill the remaining slots with the selected center
of smallest row index. Lloyd and likelihood ties use the smallest label,
empty Lloyd centers stay fixed, and equal restart objectives retain the
first restart. Matching ties follow a fixed deterministic assignment rule.

For a training fold with $m$ rows, $d$ retained items, and $r$ retained
score coordinates, screening and preprocessing cost $O(mJ)$ and $O(md)$,
respectively. Dense exact SVD costs $O(md\min\{m,d\})$. With $R$ starts
and at most $T_{\rm Lloyd}$ iterations per start, rounding costs
$O(RT_{\rm Lloyd}mKr)$. Estimating profiles costs $O(md)$ on the retained
items and $O(mJ)$ when all-item profiles are stored. Matching costs
$O(K^2d+K^3)$ and assigning a target fold $F$ costs $O(|F|Kd)$.

\subsection{Empirical EM updates}
\label{app:em-details}
The EM comparisons use all input items. If $r_{ik}$ is the posterior
probability of class $k$ for row $i$, the M-step updates the class weight
to $n^{-1}\sum_i r_{ik}$ and sets
\[
\theta_{kj}^{\rm new}=
\Pi_{[10^{-4},1-10^{-4}]}\!\left(
\frac{\sum_i r_{ik}\Omega_{ij}X_{ij}}{\sum_i r_{ik}\Omega_{ij}}
\right).
\]

Here $\Pi$ denotes interval projection. For a zero denominator, use the global
observed item mean clipped to $[.02,.98]$, or $.02$ if the item is entirely
unobserved. Floor class weights at $10^{-12}$ and renormalize. Jeffreys
smoothing is used only for initialization. Random starts perturb the clipped
global item means by independent $N(0,.12^2)$ draws, clip to $[.01,.99]$, and
use uniform weights. Reject a proposed likelihood decrease and count the
fit as nonconverged. Otherwise, stop when
$\ell_{t+1}-\ell_t\le 10^{-6}(1+|\ell_{t+1}|)$ or at 150 iterations,
where $\ell_t$ is the observed log likelihood after $t$ updates.
Final labels maximize the posterior class probabilities.

\section{Proofs of the recovery theorems}
\label{app:technical-proofs}

\subsection{Screening and fold balance}
\label{app:spectral}
\paragraph{Fold balance.}
Serfling's bound \citep{Serfling1974} yields simultaneous fold balance
outside probability $C\eta$, with deviations
$O(\sqrt{\log(2KL/\eta)/m})$.

\paragraph{Proof of Lemma 1.}
Scalar Bernstein bounds and a union bound over \(J\) items and \(L\) folds
show that when $mh\ge C\log(4JL/u)$, the following bounds hold
simultaneously with probability at least $1-u$:
\[
S_{2h}\subseteq\widehat S_{h,f}\subseteq S_{h/2},\qquad
\sup_{j\in S_{h/2}}\left|\frac{\widehat\pi_j^{(f)}}{\pi_j}-1\right|\le\delta_m(u).
\]

The core inclusion uses a
lower multiplicative tail bound, the outer inclusion uses its upper
counterpart, and division by \(\pi_j\ge h/2\) yields the relative error.

\subsection{Proof of Theorem 1}
\label{app:residual-recovery}
Condition on the training row set and fixed labels. Let
$c_j=m^{-1}\sum_{i\in T}\theta_{z_i j}$, write $\mat H=\mat H_T$ for row
centering, and on the deterministic outer set $S_o=S_{h/2}$ define
\[
R_{\gamma,ij}=\pi_j^{-\gamma}\Omega_{ij}(X_{ij}-c_j),\qquad
M_{\gamma,ij}=\pi_j^{-\gamma}(\Omega_{ij}-\pi_j).
\]

Since $c_j$ is deterministic under this conditioning, independence of the
mask and response implies the exact identity
\begin{equation}
\operatorname{Var}(R_{\gamma,ij}\mid z_i=k)
=\pi_j^{1-2\gamma}\left\{\theta_{kj}(1-\theta_{kj})
 +(1-\pi_j)(\theta_{kj}-c_j)^2\right\}.
\label{eq:residual-variance}
\end{equation}

Both $\theta_{kj}$ and $c_j$ are at most $q_+\rho$, and $\rho$ is bounded.
Thus the entries of $\mat{R}_\gamma-E\mat{R}_\gamma$ are independent, with
variance at most $C\rho\pi_j^{1-2\gamma}$ and magnitude at most $Ch^{-\gamma}$.
The rectangular
matrix variance parameter is at most $C\rho v_\gamma$.
Matrix Bernstein \citep{Tropp2012}, followed by $\|\mat{H}\|_{\rm op}=1$,
therefore bounds $\|\mat{H}(\mat{R}_\gamma-E\mat{R}_\gamma)\|_{\rm op}$ by $B_\gamma$.
The same argument without $\rho$ bounds $\|\mat{H}\mat{M}_\gamma\|_{\rm op}$ by
$B_\gamma^{\rm mask}$. Both events hold with probability at least $1-Cu$.

For each outer-set column,
\[
\widehat c_j-c_j=
\frac{\sum_{i\in T}\Omega_{ij}(X_{ij}-c_j)+1/2-c_j}
     {\sum_{i\in T}\Omega_{ij}+1}.
\]

$\sum_i E\{\Omega_{ij}(X_{ij}-c_j)\}=0$ and the summed variance is at
most $Cm\pi_j\rho$. Scalar Bernstein and a union bound control the
numerator uniformly; the observed count is at least $cm\pi_j$.
The symmetric Jeffreys terms contribute $O((m\pi_j)^{-1})$ to the
ratio, which is absorbed in $\ell_u/(mh)$. Hence
$\|\widehat{\boldsymbol c}-\boldsymbol c\|_\infty\le a_c$. Lemma~1 simultaneously implies
$S_c\subseteq\widehat S_h\subseteq S_o$ and relative propensity
error at most $\delta$. These events are intersected, without asserting
independence of the estimated center, mask or retained set.

Embed discarded columns as zeros and set
\[
D_{\gamma,jj}=\mathbf1\{\widehat\pi_j\ge h\}
 \left(\frac{\pi_j}{\max(\widehat\pi_j,h)}\right)^\gamma.
\]

With this choice of weights, the fitted residual matrix decomposes as
\begin{align*}
\widehat{\mat{R}}_\gamma
 &=\mat{P}_\gamma+\mat{H}(\mat{R}_\gamma-E\mat{R}_\gamma)\mat{D}_\gamma
   +\mat{H}\mat{M}_\gamma \mat{D}_\gamma\operatorname{diag}(\boldsymbol c-\widehat{\boldsymbol c}),\\
\mat{P}_\gamma&=\mat{H}(E\mat{R}_\gamma)\mat{D}_\gamma.
\end{align*}

The constant mean of $\mat{\Omega}\operatorname{diag}(\boldsymbol\pi^{-\gamma})$
is annihilated by $\mat{H}$. Since $\|\mat{D}_\gamma\|_{\rm op}\le(1-\delta)^{-\gamma}$,
the perturbation is at most $b_\gamma$. At $\gamma=0$, $\mat{D}_0$ is simply
the random retention projection; the same identity still handles
data-dependent selection and centering.

The signal $\mat{P}_\gamma$ has rank at most $r=K-1$. If $\widehat{\mat{P}}_\gamma$
is the rank-$r$ truncated SVD of $\widehat{\mat{R}}_\gamma$, best approximation
implies $\|\widehat{\mat{P}}_\gamma-\mat{P}_\gamma\|_{\rm op}\le2b_\gamma$ and
\[
\|\widehat{\mat{P}}_\gamma-\mat{P}_\gamma\|_F^2\le8r b_\gamma^2.
\]

Rows of $\mat{U}\mat{\Sigma}$ and $\widehat{\mat{P}}_\gamma$ have identical distances.
Let $\Delta_\gamma$ be the minimum distance between distinct random
signal centers. Then
\[
\Delta_\gamma^2\ge(1+\delta)^{-2\gamma}
\min_{k\ne\ell}\sum_{j\in S_c}\pi_j^{2-2\gamma}
                         (\theta_{kj}-\theta_{\ell j})^2
\ge c\Delta_{\gamma,c}^2.
\]

This uses only the retained core, and needs no positive $r$th singular value.
Zero padding covers fewer than $r$ retained columns. On this event,
$\Delta_{\gamma,c}>0$ makes $S_c$ nonempty: empty-support failure, with
extended loss one, is confined to the exceptional probability $Cu$.

Conditional on the realized scores, the expected D-squared initial
objective is at most $a_K=8(\log K+2)$ times the minimum within-cluster sum of squares $Q_*$
\citep[Theorem~1.1]{ArthurVassilvitskii2007}. For $Q_*>0$, Markov's
inequality bounds the probability that one start exceeds $\lambda a_KQ_*$
by $\lambda^{-1}$. Conditional independence of $R$ starts makes the
failure probability of their minimum at most $\lambda^{-R}$;
Lloyd updates cannot increase the objective. If the optimum is zero,
D-squared seeding covers every distinct score row before its sampling
mass vanishes; subsequent duplicate centers and smallest-index ties keep
the objective zero. This includes fewer than $K$ distinct rows and
requires no division by a zero optimum. For its fitted centroid
matrix $\widetilde{\mat{P}}$, comparison with the true row partition implies
$\|\widetilde{\mat{P}}-\widehat{\mat{P}}_\gamma\|_F^2
\le \lambda a_K\|\mat{P}_\gamma-\widehat{\mat{P}}_\gamma\|_F^2$.
Thus $\|\widetilde{\mat{P}}-\mat{P}_\gamma\|_F^2\le C \lambda a_K r b_\gamma^2$.
Rows whose centroid is at least $\Delta_\gamma/3$ from their true
center number at most $C \lambda a_K r b_\gamma^2/\Delta_\gamma^2$.
If this count is smaller than the smallest class, every true class
has a good row. Two good rows from distinct true classes cannot share
a fitted centroid; since there are $K$ centroids this identifies a permutation.
All good rows are then correctly assigned. In the other case the trivial
bound one is absorbed by a constant depending on class balance. This proves
the theorem, using the low-rank/center-separation route of
\citet[Lemma~4.2]{LofflerZhangZhou2021} and
\citet[Proposition~3.1]{ZhangZhou2024}.

\subsection{Proof of Theorem 2}
\label{app:refinement}
\paragraph{Profile and prior error.}
On the training-accuracy event in~\eqref{eq:initialization}, fix one fold
and align its estimated labels by the loss-minimizing permutation $\sigma_f$.
Suppress the fold superscript on profile estimates; all profile bounds in
this paragraph restrict item indices to $S_{h/2}$. Let $C_k$ and
$\widehat C_k$ be the true and aligned estimated training groups.
For $C_k=\{i\in T_f:z_i=k\}$, define observed counts, success counts
and centered success sums by \(N^0_{kj}=\sum_{i\in C_k}\Omega_{ij},\quad T^0_{kj}=\sum_{i\in C_k}\Omega_{ij}X_{ij},\quad U^0_{kj}=T^0_{kj}-\theta_{kj}N^0_{kj}.\) Fold balance and multiplicative concentration imply
\(N^0_{kj}\ge cm_f\pi_j\) on \(S_{h/2}\). Bernstein's inequality and a union bound,
with \(x_{\mathrm{prof}}=\log(8KJL/\eta)\), imply \(|U^0_{kj}|\le C\{\sqrt{m_f\pi_j\rho x_{\mathrm{prof}}} +x_{\mathrm{prof}}\}\) uniformly in $k,j,f$. Define $N_{kj}$ and $T_{kj}$ by the same sums over
$\widehat C_k=\{i\in T_f:\widehat z_i^{(f)}=k\}$. Changing at most \(m_f\epsilon_m\) memberships
implies $|N_{kj}-N^0_{kj}|\le m_f\epsilon_m$ and
$|T_{kj}-T^0_{kj}|\le m_f\epsilon_m$. This deterministic contamination
bound does not treat the adaptive error set as independent. Write the true
count bound as $N^0_{kj}\ge c_1m_f\pi_j$. Since $\rho$ is bounded,
choosing the constant in $r_m\le c$ small enough ensures
$\epsilon_m\le c_1\pi_j/2$ on $S_{h/2}$, and hence
\[
N_{kj}\ge N^0_{kj}-m_f\epsilon_m\ge(c_1/2)m_f\pi_j.
\]

For the Jeffreys estimate the exact identity is
\[
\widehat\theta_{kj}-\theta_{kj}=
\frac{U^0_{kj}+(T_{kj}-T^0_{kj})
-\theta_{kj}(N_{kj}-N^0_{kj})+1/2-\theta_{kj}}{N_{kj}+1}.
\]

Applying the count bounds to the numerator and denominator gives a uniform
bound on the relative profile error:
\[
\max_{k,j,f}\frac{|\widehat\theta_{kj}-\theta_{kj}|}{\rho}
\le C\!\left[
\frac{\epsilon_m}{h\rho}
+\sqrt{\frac{x_{\mathrm{prof}}}{mh\rho}}
+\frac{x_{\mathrm{prof}}}{mh\rho}\right]
=Cr_m.
\]

Define the true training fraction
$\alpha_{k,f}=m_f^{-1}\sum_{i\in T_f}\mathbf1\{z_i=k\}$;
the estimated labels use the same relabeling by $\sigma_f$.
Direct algebra for the Jeffreys class weights yields
\(\max_{k,f}|\widehat\alpha_k^{(f)}-\alpha_{k,f}|\le Ca_m,\qquad
\max_{k,f}|\log\widehat\alpha_k^{(f)}-\log\alpha_{k,f}| \le C a_m.\)
\paragraph{Cross-fold profile matching.}
The relative profile event and regular-profile bounds imply, on every common
empirical subset of \(S_{h/2}\), \(d_S^2(\widehat{\boldsymbol{\theta}}_k^{(f)},\boldsymbol{\theta}_k) \le C|S_{h/2}|\rho r_m^2=e_m.\) Any two retained sets contain \(S_{2h}\). The triangle inequality in the
square-root embedding implies that a correct pair has distance at most
$2\sqrt{e_m}$; every incorrect pair has distance at least
$\sqrt{D^{\mathrm{prof}}_{*,2h}}-2\sqrt{e_m}$. Under condition~\eqref{eq:matching} the latter
strictly exceeds the former. Thus the correct column is the strict minimum
in every row of the squared-distance cost matrix. These columns are
distinct, so they form the unique minimum-cost assignment and align all
folds to one global permutation, without a factor $K$ in~\eqref{eq:matching}.

\paragraph{Held-out score error.}
For $i\in F_f$, define the oracle score on the fold-balance event by
$L_{ik}^{(f)}=\log\alpha_{k,f}+\sum_{j\in\widehat S_{h,f}}\Omega_{ij}
\{X_{ij}\log\theta_{kj}+(1-X_{ij})\log(1-\theta_{kj})\}$,
and assign it arbitrary values off that event.
The estimated score $\widehat L_{ik}^{(f)}$ replaces these profiles and
proportions with their Jeffreys estimates, analytically relabeled by $\sigma_f$.
Write $f(i)$ for the unique target fold containing row $i$.
Let \(\mathcal G_f\) be the intersection of the training-side fold-balance,
propensity, initializer, profile, and prior events. It is measurable with
respect to the complementary training data, the split, and training
randomization; it excludes the target row. The profile event implies \(|\log\widehat\theta_{kj}-\log\theta_{kj}|\le Cr_m,\qquad |\log(1-\widehat\theta_{kj})-\log(1-\theta_{kj})| \le C\rho r_m.\) Conditional Bernstein bounds for the number \(P_i^+\) of observed positives
and \(N_i^{\mathrm{obs}}\) of observed cells on the deterministic outer set $S_{h/2}$ imply, uniformly over
held-out rows, \(P_i^++\rho N_i^{\mathrm{obs}} \le C\{D_{h/2}+\sqrt{D_{h/2}y_n}+y_n\}.\) Let $E_{\mathrm{row}}$ denote the simultaneous event that this row-count
bound holds for every held-out row, with counts taken on the deterministic
outer set $S_{h/2}$. Its complement has probability at most $C\eta$ by a
union bound. Combining this event with the class-weight error yields
\(\max_{i,k}|\widehat L_{ik}^{(f(i))}-L_{ik}^{(f(i))}| \le R_m^J.\)

\paragraph{Probability of a classification error.}
Let $\mathcal T_f$ be the training sigma-field defined alongside
condition~\eqref{eq:initialization}.
Let $E_i$ be the row-count event and $\mathcal N$ the simultaneous
nuisance and correct-matching event. On $\mathcal N$ no failed return
occurs. If the final algorithm misclassifies a row after the common
global relabeling, some incorrect class has score at least as large as
the true class. This implication includes score ties and does not require
argmax tie-breaking to commute with a permutation. Consequently,
\begin{align*}
P\{\operatorname{err}(\widehat{\boldsymbol z},\boldsymbol z)>0\}
&\le P(\mathcal N^c)+P(E_{\mathrm{row}}^c)\\
&\quad+\sum_f\sum_{i\in F_f}\sum_{\ell\ne z_i}
P\{\widehat L_{i\ell}^{(f)}\ge\widehat L_{iz_i}^{(f)},\mathcal G_f,E_i\}.
\end{align*}

For a held-out subject from class \(k\), let
\(\Lambda_{k\ell}^{(f)}=L_{ik}^{(f)}-L_{i\ell}^{(f)}\), including the true
fold class weights. Target-fold independence implies, for \(u\in[0,1]\), \(E_k\{\exp(-u\Lambda_{k\ell}^{(f)})\mid\mathcal T_f\} =\left(\frac{\alpha_{\ell,f}}{\alpha_{k,f}}\right)^u \prod_{j\in\widehat S_{h,f}}c_{k\ell j}(1-u).\) The score comparison $\widehat L_{i\ell}^{(f)}\ge\widehat L_{ik}^{(f)}$ implies
$\Lambda_{k\ell}^{(f)}\le2R_m^J$ on $\mathcal G_f\cap E_i$.
Before taking the training-conditional exponential moment, use the event inclusion
\[
\{\widehat L_{i\ell}^{(f)}\ge\widehat L_{ik}^{(f)},\mathcal G_f,E_i\}
\subseteq\{\Lambda_{k\ell}^{(f)}\le2R_m^J\}\cap\mathcal G_f.
\]

No independence
between $E_i$ and the likelihood ratio is needed. Because
\(S_{2h}\subseteq\widehat S_{h,f}\) and each coefficient is at most one,
Markov's inequality yields
\[
P_k\{\widehat L_{i\ell}^{(f)}\ge\widehat L_{ik}^{(f)},\mathcal G_f,E_i\}
\le C\exp\!\left[
\sum_{j\in S_{2h}}\log c_{k\ell j}(1-u)+2uR_m^J\right].
\]

Optimizing in \(u\) yields the penalized Chernoff exponent. Taking \(u=1/2\)
yields the R\'enyi-based bound. If \(u_*\) optimizes the oracle coefficient,
evaluating the penalized supremum at \(u_*\) yields
\(\mathcal C_{k\ell,S}(\xi)\ge C^{\mathrm{obs}}_{k\ell,S}-2\xi\).
Union bounds over rows, competitors, folds, and nuisance-event complements
complete the theorem. Matching failure is charged once through
$P(\mathcal N^c)$. Only the target fold's training event enters its
rowwise Chernoff step; the held-out response remains random in this step.

\subsection{Estimates for the proof of Theorem 3}
\label{app:sharp-outline}
Condition on a training fold of size $m\asymp n$ and its fixed labels.
Write $\mat V=\mat\Omega\odot\mat X$, $P_{ij}=\pi_j\theta_{z_i j}$,
$\mat E=\mat V-\mat P$, and
$\mat M=\mat\Omega-\mathbf1\boldsymbol\pi^\top$, where
$\boldsymbol\pi=(\pi_j)_{j=1}^J$. On the training event established in
Online Resource~1, Section~\ref*{app:sharp-initialization},
\[
\|\mat E\|\le C\sqrt{d_n},\quad \|\mat M\|\le C\sqrt{n\pi},\quad
\max_j\sum_iV_{ij}\le Cd_n,\quad \max_j\sum_i\Omega_{ij}\le Cn\pi.
\]

On the same foldwise events, the retained sets omit at most $C/\pi$
items per fold and the training loss is $O(d_n^{-1/2})$. Taking the
concentration failure exponent to be at least two and summing the
restart bounds over folds gives total exceptional probability at most
$Cn^{-2}+CLd_n^{-R_0/2}$.

\paragraph{Aggregate profile error and inverse bounds.}
On the matrix and count event just stated, consider any partition with
permutation loss at most $\epsilon=o(1)$, analytically aligned to the
truth. For group $k$ write its membership indicator as $w_i=\mathbf1\{\widehat z_i=k\}$,
its size as $\widehat n_k=\sum_iw_i\asymp n$, and define the observed and
success counts $N_{kj},T_{kj}$, Jeffreys estimate $a_{kj}=\widehat\theta_{kj}$,
and good-denominator item set $G_k$ by
\[
N_{kj}=\boldsymbol w^\top\mat\Omega_{\cdot j},\quad
T_{kj}=\boldsymbol w^\top\mat V_{\cdot j},\quad
a_{kj}=\frac{T_{kj}+1/2}{N_{kj}+1},\quad
G_k=\{j:N_{kj}\ge\widehat n_k\pi_j/2\}.
\]

Put $r_{\rm agg}^2=\epsilon^2+d_n^{-1}$ for the aggregate-error budget.
Use the single-item distance
$d_{\rm Ber}^2(a,p)=(\sqrt a-\sqrt p)^2+(\sqrt{1-a}-\sqrt{1-p})^2$,
as in Section~\ref{sec:method}. The following bounds hold simultaneously
for all these adaptive partitions:
\begin{align}
|G_k^c|&\le C/\pi,&
\sum_{j\in G_k}(a_{kj}-\theta_{kj})^2&\le CJ\rho^2r_{\rm agg}^2,\notag\\
\sum_j d_{\rm Ber}^2(a_{kj},\theta_{kj})&\le CJ\rho r_{\rm agg}^2,&
a_{kj}&\le C\rho\quad(j\in G_k),\notag\\
\rho/a_{kj}&\le Cd_n,&(1-a_{kj})^{-1}&\le Cd_n\quad(\text{all }j).
\label{eq:sharp-aggregate}
\end{align}

To prove this, define
$\bar\theta_{kj}=\widehat n_k^{-1}\sum_iw_i\theta_{z_i j}$.
At most $m\epsilon$ foreign rows enter group $k$, so
$|\bar\theta_{kj}-\theta_{kj}|\le C\epsilon\rho$.
Let $\boldsymbol F_k=\boldsymbol w^\top\mat E$ and
$\boldsymbol D_k=\boldsymbol w^\top\mat M$ be the grouped response and mask
fluctuations, with entries $F_{kj}$ and $D_{kj}$. Pathwise norm bounds imply
\[
\|\boldsymbol F_k\|_2^2\le Cnd_n,\qquad
\|\boldsymbol D_k\|_2^2\le Cn^2\pi.
\]

Since $N_{kj}=\widehat n_k\pi_j+D_{kj}$, every bad denominator costs
at least $cn^2\pi^2$ in the second sum; hence $|G_k^c|\le C/\pi$.
The profile error satisfies the exact identity
\[
a_{kj}-\theta_{kj}
=\frac{\widehat n_k\pi_j(\bar\theta_{kj}-\theta_{kj})
       +F_{kj}-\theta_{kj}D_{kj}+1/2-\theta_{kj}}{N_{kj}+1}.
\]

Combining this identity with $N_{kj}\ge cn\pi$ on $G_k$ gives
\[
\sum_{j\in G_k}(a_{kj}-\theta_{kj})^2
\le C\left[J\epsilon^2\rho^2+
 \frac{nd_n+\rho^2n^2\pi+J}{n^2\pi^2}\right]
\le CJ\rho^2r_{\rm agg}^2.
\]

The success-count bound also implies $a_{kj}\le C\rho$ on $G_k$.
For regular $p$, $d_{\rm Ber}^2(a,p)\le(a-p)^2\{p^{-1}+(1-p)^{-1}\}
\le C(a-p)^2/\rho$. Bad denominators contribute at most
$2|G_k^c|\le C/\pi=O(J\rho/d_n)$, proving the Hellinger bound.
Finally, smoothing and the two total-count bounds imply
\[
\frac\rho{a_{kj}}\le2\rho(N_{kj}+1)\le Cd_n,\qquad
\frac1{1-a_{kj}}=1+\frac{T_{kj}+1/2}{N_{kj}-T_{kj}+1/2}
\le2(T_{kj}+1)\le Cd_n.
\]

The success count, rather than the observation count, controls the second
inverse. These identities include zero observations and all-success groups.

\paragraph{Fractional moments and likelihood reversals.}
Under the bounds in Equation~\eqref{eq:sharp-aggregate} and
$\epsilon\le C_0d_n^{-1/2}$, the full-product oracle order $u_*$ lies in
$[\beta,1-\beta]$ for a fixed $\beta>0$. For every item set $S$ and every
$u$ in this interval, the plug-in coefficients satisfy
\[
\sum_{j\in S}\log\widehat c_j(u)
\le \sum_{j\in S}\log c_{k\ell j}(1-u)+C(\sqrt{d_n}+d_n^u).
\]

The coefficient for item $j$ is
$\widehat c_j(u)=1-\pi_j+\pi_j\{p_j(b_j/a_j)^u+(1-p_j)[(1-b_j)/(1-a_j)]^u\}$,
where $p_j=\theta_{kj}$, $a_j=\widehat\theta_{kj}$ and $b_j=\widehat\theta_{\ell j}$.

The oracle coefficients are bounded below, so each positive log-coefficient
difference is at most $C\pi_j$ times the positive difference in the response
terms. First group items with a bad denominator, then those with an estimate
close to one, and then those with an estimate below half its true value,
assigning each item to the first applicable set. The count and squared-error
bounds give respective set sizes at most $C/\pi$, $CJ\rho^2r_{\rm agg}^2$,
and $CJr_{\rm agg}^2$. The inverse bounds give per-item excess at most
$C\pi_jd_n^u$ in the first two sets and $C\pi_j\rho d_n^u$ in the third,
where the failure term is Lipschitz. Each set contributes at most $Cd_n^u$
because $r_{\rm agg}^2\le C/d_n$. On the remaining items, the estimates are
of order $\rho$ and bounded away from one, so the response term is Lipschitz
in $a_j,b_j$. The squared-error bound and Cauchy--Schwarz give a total excess
of at most $C\pi J\rho r_{\rm agg}\le C\sqrt{d_n}$ on these items. Since
$u_*\le 1-\beta$, the combined excess at $u_*$ is $o(d_n)$ as
$d_n\to\infty$. Online Resource~1, Section~\ref*{app:sharp-moments}, gives
the full derivation.

For $d_n\to\infty$, a hardest-pair full-row log likelihood ratio $\Lambda$
satisfies
\[
P_k\{\Lambda<0\},\ P_\ell\{\Lambda>0\}
\ge c\exp\{-C_*^{\rm obs}-C\sqrt{d_n}\}.
\]

For bounded $d_n$, the success probability under the disjoint-pair prior
is at most $K!(1-c)^M$ up to permutation, where $M\asymp n$ is the number of
pairs and $c>0$ bounds each pair's Bayes error from below.
Online Resource~1, Section~\ref*{app:sharp-converse}, proves both estimates.

\paragraph{Uniform and aggregate error scales.}
\label{app:route-comparison}
Under comparable observation probabilities, put $s_n=\log n/(n\pi\rho)$
and $t_n=\log n/(n\pi^2\rho^2)$. If $t_n\to0$, Theorem~1 gives
$\epsilon_m\lesssim\lambda_n(s_n+s_n^2)$. Taking
$\lambda_n=t_n^{-1/2}$ makes
$\epsilon_m/(\pi\rho)=O(\sqrt{t_n}(1+s_n))=o(1)$ and the restart failure
probability vanish, so the uniform-profile argument of Theorem~2 applies
under its matching and concentration conditions. At the sharp threshold,
with $\pi_j=1$ and $\rho\asymp\log n/n$, this sufficient scale is unavailable,
while the aggregate argument above controls the exponent loss.

\section{Proofs of the recovery corollaries}
\label{app:corollary-proofs}

\subsection{Proof of Corollary 1}
\label{app:common-threshold-proof}
For each fixed $\gamma$, the upper bound in Theorem~3 tends to zero above
$(1+\delta)\log n$ because $C_*^{\rm obs}\asymp d_n$ and its exponent
correction is $o(d_n)$. The converse applies to each initialization
representation. For the column-centered raw input, the following
rescaling transfers the same initialization bounds.

On $S_{h/2}$, define the oracle-weighted raw matrix $\mat{Y}^{-,f}\in\mathbb R^{m_f\times|S_{h/2}|}$
on training rows $T_f$, with oracle mean
$\mat{P}^{-,f}=E(\mat{Y}^{-,f}\mid\boldsymbol z,T_f)$.
The superscript $-$ marks oracle weighting. Define the diagonal rescaling
matrix $\mat D_f$ (off-diagonal entries zero) and empirical input
$\widehat{\mat Y}^{(f)}$ by
\begin{align*}
Y^{-,f}_{ij}&=\frac{\Omega_{ij}X_{ij}}{\sqrt{\pi_j}},\\
D_{f,jj}&=\mathbf{1}\{\widehat\pi_j^{(f)}\geq h\}
\sqrt{\frac{\pi_j}{\max\{\widehat\pi_j^{(f)},h\}}},\\
\widehat{\mat{Y}}^{(f)}&=\mat{Y}^{-,f}\mat{D}_f,\\
\mat{H}_{T_f}(\widehat{\mat{Y}}^{(f)}-\mat{P}^{-,f}\mat{D}_f)
&=\mat{H}_{T_f}(\mat{Y}^{-,f}-\mat{P}^{-,f})\mat{D}_f.
\end{align*}

Since $h=\pi/2$ and $\pi_j\le C_\pi\pi$, the definition gives
$\|\mat D_f\|\le\sqrt{2C_\pi}$ for every realization. On the training event
of Online Resource~1, Section~\ref*{app:sharp-initialization}, the centered
success matrix $\mat E$ of Appendix~\ref{app:sharp-outline} satisfies
$\|\mat E\|\le C\sqrt{d_n}$. Hence
$\|\mat Y^{-,f}-\mat P^{-,f}\|\le\pi^{-1/2}\|\mat E\|$, and the identity
above bounds the centered Raw perturbation by $C\sqrt{d_n/\pi}$.

On the same event, $|\widehat S_f^c|\le C/\pi=o(J)$ and
$\widehat\pi_j^{(f)}\le C\pi$ on $\widehat S_f$. The scaled contrasts
$q_{kj}-q_{\ell j}$ are bounded, so deleting these $o(J)$ items removes
only $o(J)$ from each pair's sum of squared contrasts. By the separation
assumption of Theorem~3 and $\pi_j\ge\pi$, the minimum squared distance
between retained Raw class centers after centering satisfies
\[
\min_{k\ne\ell}\sum_{j\in\widehat S_f}
  \frac{\pi_j^2}{\widehat\pi_j^{(f)}}(\theta_{kj}-\theta_{\ell j})^2
  \ge cJ\pi\rho^2.
\]

These perturbation and separation bounds give the $O(d_n^{-1/2})$ training
loss bound with the same restart failure probability as in that section.
The profile estimation, fold alignment, and likelihood refinement steps of
Theorem~3 then apply to the Raw partition.

\subsection{Proof of Corollary 2}
\label{app:rare-tail-proof}
Use the core/tail notation $S_c,J_c,T_n,p_n$ of Corollary~2.
Let $B_{\rm tail}$ be the event that at least one training fold retains a tail item.
Since tail propensities are at most $p_n/8$ and $h=p_n/2$, a binomial
upper-tail bound and a union bound give $P(B_{\rm tail})\le LT_n\exp(-cmp_n)$.
Write $\operatorname{err}_{\rm full}$ and $\operatorname{err}_{\rm core}$
for the full and core-only permutation losses, with loss one for a return of $\bot$.
Couple the full and core-only algorithms using the same core data, folds
and randomization. Their outputs agree on $B_{\rm tail}^c$, so
\[
P(\operatorname{err}_{\rm full}>0)\le
P(B_{\rm tail})+P(\operatorname{err}_{\rm core}>0).
\]

On the core, $C_{*,S_c}^{\rm obs}\asymp np_n\rho\ge c\log n$, so
Corollary~1 applies with dimension $J_c\asymp n$ and controls the second term.
Bounded $\rho$ gives $np_n\to\infty$; $m\asymp n$ and
$\log(T_n+1)=o(np_n)$ make the first term vanish. This proves Corollary~2
without conditioning on tail exclusion or using discarded tail information.

\subsection{Proof of Corollary 3}
\label{app:random-retention-proof}
Write $S_o=S_{h/2}$ and $S_+=S_{2h}$, so $S_c\subseteq S_+\subseteq S_o$
and $|S_o|\le J\asymp n$. On $S_o$ the propensities are between $h/2$
and $C_\pi h$. Put $s_n=\log n/(nh\rho)$ and
$t_n=\log n/(nh^2\rho^2)\to0$. Theorem~1 implies, uniformly over folds,
\[
v_\gamma\lesssim nh^{1-2\gamma},\qquad
\Delta_{\gamma,c}^2\gtrsim n\rho^2h^{2-2\gamma},\qquad
\frac{b_\gamma^2}{m\Delta_{\gamma,c}^2}=O(s_n+s_n^2).
\]

In this display the theorem's
$\Delta_{\gamma,c}$ uses $S_+$; the lower bound uses only $S_c$.
The center term is again $o(B_\gamma)$ because $mh\rho/\log n\to\infty$.
Choose $u_n=1/(L\log n)$, $\lambda_n=t_n^{-1/2}$ and
$\epsilon_m=C\lambda_n(s_n+s_n^2)$. Fold balance and initialization
hold outside probability
$C/\log n+Lt_n^{R_0/2}$, while $\epsilon_m/(h\rho)=o(1)$.
With $\eta_n=C/\log n+L\sqrt{t_n}$, the nuisance scale satisfies
$r_m=O(\sqrt{t_n}+\sqrt{s_n}+s_n)=o(1)$.

The fixed core supplies
$D^{\mathrm{prof}}_{*,2h}\ge c n\rho$, whereas
$e_m\le C J\rho r_m^2=o(n\rho)$, proving condition~\eqref{eq:matching}.
Regular scaling and the observation-probability bounds imply
\[
I^{\mathrm{obs}}_{*,S_+}\asymp nh\rho,\qquad
D_{h/2}=\rho\sum_{j\in S_o}\pi_j\asymp nh\rho\gg\log n.
\]

Consequently $R_m^J/I^{\mathrm{obs}}_{*,S_+}=O(r_m)+o(1)=o(1)$,
and Theorem~2 yields exact recovery. Its training-side events allow the retained set to vary between the core
and the outer set.
For an item with $\pi_j=h$, its training count has binomial distribution
$\operatorname{Bin}(m_f,h)$, with $m_f$ trials, success probability $h$ and diverging variance. The rounding of
$m_fh$ changes the standardized threshold by $o(1)$, proving the
stated limiting selection probability.

\FloatBarrier

\backmatter
\section*{Statements and Declarations}

\paragraph*{Funding.}
This work was supported by the National Natural Science Foundation of China
(Grant Nos. 12301358 and 42450275) and the Fundamental Research Funds for the
Central Universities at Wuhan University (Grant No. 2042025kf0051).

\paragraph*{Competing interests.}
The authors declare no competing interests.

\paragraph*{Data availability.}
The PISA data are available from the
\href{https://www.oecd.org/en/data/datasets/pisa-2022-database.html}{OECD PISA 2022 database}.
The simulation data were generated according to the designs described in the article.

\paragraph*{Use of generative AI.}
Generative AI tools (including OpenAI Codex) were used solely to assist with simulation coding and language polishing. All generated code and text were verified and edited by the authors, who take full responsibility for the integrity and accuracy of the work.

\bibliography{references}

\end{document}